\documentclass[12pt]{article}
\usepackage[utf8]{inputenc}
\usepackage{amsmath}
\usepackage{amsfonts}
\usepackage{amssymb}
\usepackage{amsthm}
\usepackage{relsize}
\usepackage{array}
\usepackage{multirow}
\usepackage{tabu}
\usepackage{soul}
\usepackage{graphicx}
\usepackage{epigraph}
\usepackage{float}
\usepackage[english]{babel}
\usepackage[usenames, dvipsnames]{color}
\usepackage{fancyhdr}
\usepackage[Sonny]{fncychap}
\usepackage{tikz-cd}
\usepackage{tkz-euclide}
\usepackage{sseq}
\usepackage{newfloat}
\usepackage{amsmath,mathtools}
\usepackage{pgfplots}
\usepackage{geometry}
 \usepackage{sectsty}
\sectionfont{\fontsize{14}{15}\selectfont}

\usepackage[all]{xy}
\DeclareFloatingEnvironment[fileext=dia]{diagram}
\theoremstyle{definition}

\newtheorem{remark}{Remark}

\usepackage[most]{tcolorbox}

\makeatletter
\newtcolorbox{note}[1][]{%
	breakable,
	enhanced jigsaw, % better frame drawing
	borderline west={3pt}{0pt}{black!10!white}, % straigh vertical line at the left edge
	borderline south={1pt}{0pt}{black!10!white}, 
	borderline east={1pt}{0pt}{black!10!white},
	borderline north={1pt}{0pt}{black!10!white},
	sharp corners, % No rounded corners
	boxrule=0pt, % no real frame,
	attach title to upper, % Move the title into the box,
	left=0pt,
	right=0pt,
	top=0pt,
	bottom=0pt,
	boxsep=5pt,
	colback=white,
	frame hidden,
	#1
}
\makeatother

\makeatletter
\newtcolorbox{note1}[1][]{%
	breakable,
	enhanced jigsaw, % better frame drawing
	sharp corners, % No rounded corners
	boxrule=0pt, % no real frame,
	attach title to upper, % Move the title into the box,
	fontupper=\linespread{1.1}\fontfamily{qpl}\selectfont,
	fontlower=\linespread{1.1}\fontfamily{qpl}\selectfont, 
	left=0pt,
	right=0pt,
	top=0pt,
	bottom=0pt,
	boxsep=3pt,
	colback=green!3!white,
	frame hidden,
	before skip=10pt plus 2pt,after skip=10pt plus 2pt,
	#1
}
\makeatother

\makeatletter
\newcommand\tabfill[1]{%
	\dimen@\linewidth
	\advance\dimen@\@totalleftmargin
	\advance\dimen@-\dimen\@curtab
	\parbox[t]\dimen@{#1\ifhmode\strut\fi}%
}
\makeatother

\usepackage{tabularx}
\usepackage{imakeidx}
\usepackage{hyperref}
\hypersetup{colorlinks={true},linkcolor={blue},citecolor={blue},urlcolor={blue}, linktoc=all}  
 
 \usepackage[capitalize, nameinlink]{cleveref}
 \crefdefaultlabelformat{#2\textbf{#1}#3}
 \crefname{figure}{Figure}{Figures} 
 \Crefname{figure}{Figure}{Figures}
 \crefname{table}{Table}{Tables}
 \Crefname{table}{Table}{Tables}
 \crefname{section}{\S\hspace{-1mm}}{\S\hspace{-1mm}}
 \Crefname{section}{\S\hspace{-1mm}}{\S\hspace{-1mm}}
 \crefname{equation}{}{}
 \Crefname{equation}{}{}
 \crefname{example}{Geometric Pattern}{Geometric Patterns} 
 \Crefname{example}{Geometric Pattern}{Geometric Patterns}

 \usepackage{graphicx,tipa}% http://ctan.org/pkg/{graphicx,tipa}

 \usepackage{footnote}

\begin{document}

\title{\textbf{Quadratic Equations in Elamite Mathematics}}

\author{Nasser Heydari\footnote{Email: nasser.heydari@mun.ca}~ and  Kazuo Muroi\footnote{Email: edubakazuo@ac.auone-net.jp}}

\maketitle

\begin{abstract}
In this article, we study  some of quadratic equations and their solutions found  in the Susa Mathematical Texts (\textbf{SMT}). We show that the Susa scribes used this group of equations in different problems and took a  standard approach, known as completing the square, to find   solutions.   
\end{abstract}

\section{Introduction}
Quadratic equations and their applications appear in several texts of the \textbf{SMT}. We  consider here  the quadratic equations  in \textbf{SMT No.\,5}, \textbf{SMT No.\,6}, \textbf{SMT No.\,20} and \textbf{SMT No.\,21}. These texts were inscribed by Elamite scribes between 1894--1595 BC on   26 clay tablets excavated from Susa in  southwest Iran by French archaeologists in 1933. The texts of all the tablets,  along with their interpretations, were first published in 1961 (see \cite{BR61}).   

\textbf{SMT No.\,5}\footnote{The reader can see   this tablet on the website of the Louvre's collection. Please see \url{https://collections.louvre.fr/ark:/53355/cl010185655} for obverse  and   reverse.}  is a collection of short problems which are intended to  teach    methods  for  solving  quadratic equations\index{quadratic equation}. The text of this tablet consists of two parts: the mathematical part and the colophon part.  According to the colophon part,  the text originally contained 262 problems, of which   184 or more  are clearly legible. This tablet was written by two Susa scribes\index{Susa scribes} named {\fontfamily{qpl}\selectfont \textit{D\={a}di-\v{s}u\v{s}nak}} (literally: a favorite of the city-god\index{city-god of Susa} of Susa\index{city-god of Susa}) and {\fontfamily{qpl}\selectfont \textit{S\^{\i}n-i\v{s}meanni}} (literally: the moon-god\index{moon-god of Susa} heard my request).

The content of  text of  \textbf{SMT No.\,6}\footnote{The reader can see  this tablet on the website of the Louvre's collection. Please see \url{https://collections.louvre.fr/ark:/53355/cl010185656} for obverse  and   reverse.}   is similar to that of \textbf{SMT No.\,5}, in that it contains  a collection of   short problems regarding quadratic equations  without   solutions. The first part deals with 39 quadratic equations  in the form of $x^2 + ax = b$, of which at least 32 problems remain. The second part of the text consists of 30 problems  in the form of  $x^2 - ax = b$, of which about 20 problems are comprehensible.

\textbf{SMT No.\,20}\footnote{The reader can see   this tablet on the website of the Louvre's collection. Please see \url{https://collections.louvre.fr/ark:/53355/cl010186430} for obverse  and   reverse.}  contains two similar algebraic problems, one on the obverse and the other on the reverse of the tablet. Both problems solve quadratic equations\index{quadratic equation} by completing the square\index{completing the square method}. They  deal  with the geometric figure {\fontfamily{qpl}\selectfont  \textit{apusamikkum}}\index{apusamikkum@\textit{apusamikkum} (geometrical figure)} that we have described elsewhere in  \textbf{SMT No.\,3} (see \cite{HM22-3}).

The text  of \textbf{SMT No.\,21}\footnote{The reader can see  this tablet on the website of the Louvre's collection. Please see \url{https://collections.louvre.fr/ark:/53355/cl010186431} for obverse  and   reverse.}   consists of two similar problems dealing with certain quadratic equations\index{quadratic equation} concerning the area\index{area of an \textit{apusamikkum}} of {\fontfamily{qpl}\selectfont \textit{apusamikkum}}\index{apusamikkum@\textit{apusamikkum} (geometrical figure)}. The first problem is comparatively well preserved, but the second, which might have started on the obverse and continued on the reverse of the tablet, is mostly lost--leaving only the last steps of the solution.

\section{Quadratic Equations}
Quadratic equations belong to the category of polynomial equations. In such equations, one is usually interested in finding the zeros of a polynomial. Geometrically, these values are the $x$-coordinates of the intersection of the graph of the polynomial with the $x$-axis. Although finding the zeros of polynomials is not possible in many cases, for   polynomials of specific degrees this can be achieved. 

\subsection{General Form and Solutions}\label{GFS}
 The general form of  quadratic equations is 
\begin{equation}\label{eq-a}
	ax^2+bx+c=0
\end{equation} 
where $a\neq 0,b,c$ are real numbers and $x$ the unknown variable. The standard way of finding the real solutions of  a quadratic equation  is to use the \textit{quadratic formula}
\begin{equation}\label{eq-b}
x=\frac{-b\pm \sqrt{b^2-4ac}}{2a},
\end{equation} 
provided that the \textit{discriminant} of the equation 
 $$\Delta = b^2-4ac$$
  is non-negative. Note the if $\Delta <0$, one can consider the complex solutions of the equation.

   The formula \cref{eq-b} can be obtained by using an approach called \textit{completing the square}. In this method, one can use the binomial identity
\begin{equation}\label{eq-c}
(A+B)^2=A^2+2AB+B^2
\end{equation}
to obtain formula \cref{eq-b}.  If we multiply both sides of  \cref{eq-a} with $a\neq 0$  and use $A=ax, B=\frac{b}{2}$, we can  write as follows: 
\begin{align*}
 &~~	a^2x^2+abx+ac   =0 \\
\Longrightarrow~~&~~ 	(ax)^2+2\left(\frac{b}{2}\right)(ax)+ac  =0 \\
\Longrightarrow~~&~~ 	(ax)^2+2\left(\frac{b}{2}\right)(ax)+\left(\frac{b}{2}\right)^2+ac =\left(\frac{b}{2}\right)^2 \\
\Longrightarrow~~&~~ 		\left(ax+\frac{b}{2}\right)^2  =\frac{b^2}{4}-ac \\
\Longrightarrow~~&~~ 			 ax+\frac{b}{2}  =\pm \sqrt{\frac{b^2}{4}-ac} \\
\Longrightarrow~~&~~ 			 ax  =- \frac{b}{2}\pm \frac{\sqrt{b^2-4ac}}{2} \\
\Longrightarrow~~&~~ 			 x  =  \frac{-b\pm \sqrt{b^2-4ac}}{2a}.
\end{align*}
In step 2, we added the term $\left(\frac{b}{2}\right)^2$ to both sides of the equation to complete the right-hand side of the binomial identity \cref{eq-b} with $A=ax$ and $B=\frac{b}{2}$. This is the reason this method is usually called completing the square.

\subsection{Quadratic System of Equations}\label{QSE}
Sometimes the sum and product of two quantities are known and we need to find the values of the quantities. For example, if the area and the perimeter of a rectangle are known, how can we determine its length and width. Such problems lead to solving  quadratic equations. Strictly speaking, let $x$ and $y$ be two unknown numbers such that
\begin{equation}\label{eq-ca}
\begin{cases}
	x+y=p\\
	xy=q
\end{cases}
\end{equation}
 where $p,q$ are known numbers. If we substitute the value of $y=\frac{q}{x}$ in the first equation and simplify, we get the following quadratic equation with respect to $x$:
 \begin{equation}\label{eq-cb}
 x^2-px+q=0.
 \end{equation} 
By the quadratic formula, one can solve this equation to find the values of $x,y$:
\[ x=\frac{p+\sqrt{p^2-4q}}{2}~~~\text{and}~~~y=\frac{p-\sqrt{p^2-4q}}{2} \]
provided that $p^2-4q\geq 0$.

The  system of equations \cref{eq-ca}  belongs to  quadratic systems of equations because its solution requires  solving a quadratic equation.   

Another way to solve \cref{eq-ca} is to use the algebraic identity  
 \begin{equation}\label{eq-ce}
(A+B)^2-(A-B)^2=4AB.
\end{equation}  
If we set $A=\frac{x}{2}$ and  $B=\frac{y}{2}$,  it follows from   \cref{eq-ca}   that 
\begin{equation}\label{eq-cda}
xy=\left(\frac{x+y}{2}\right)^2-\left(\frac{x-y}{2}\right)^2,
\end{equation}
which consequently  gives us the following two equations
\begin{equation}\label{eq-cdb}
	\frac{x+y}{2} =\sqrt{xy+\left(\frac{x-y}{2}\right)^2},
\end{equation}
and
\begin{equation}\label{eq-cdc}
	\frac{x-y}{2} =\sqrt{\left(\frac{x+y}{2}\right)^2-xy}.
\end{equation}
Since 
$$x=\frac{x+y}{2}+\frac{x-y}{2}$$
and
$$y=\frac{x+y}{2}-\frac{x-y}{2},$$
 it follows from \cref{eq-ca}, \cref{eq-cdb} and \cref{eq-cdc} that
\[ x=\ \frac{p+\sqrt{p^2-4q}}{2}~~~~~\text{and}~~~~~y=  \frac{p-\sqrt{p^2-4q}}{2},\] 
or
\[ x=\ \frac{p-\sqrt{p^2-4q}}{2}~~~~~\text{and}~~~~~y=  \frac{p+\sqrt{p^2-4q}}{2}.\] 

\subsection{Geometric Viewpoint}
In general, any function like $f(x)=ax^2+bx+c$ where $a\neq0,b,c$ are real numbers, is a quadratic function. The graph of such a function is a parabola in the $xy$-coordinate system. It can be shown that the graph of any quadratic function is obtained from the graph of the \textit{unit parabola} $y=x^2$. This can be proved using completing the square as above. In fact,
\[ f(x)=a\left(x+\frac{b}{2a}\right)^2+\frac{4ac-b^2}{4a} \]
This formula shows that the graph of $f(x)$ can be formed by translating, reflecting, compressing or stretching the  graph of $y=x^2$. 

\begin{figure}[H]
	\centering
	\includegraphics[scale=1]{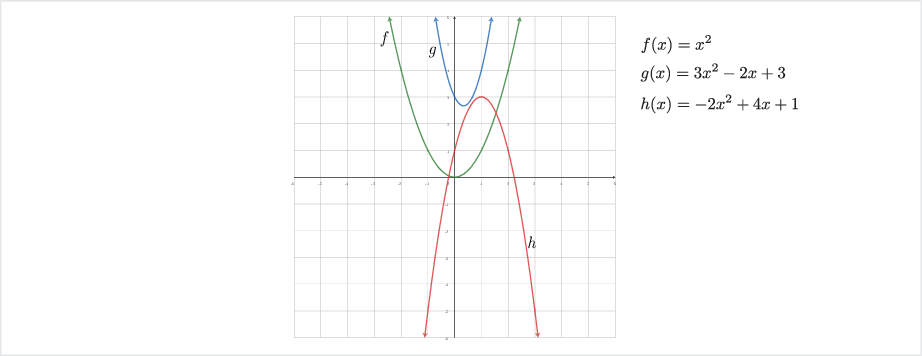}
	%\vspace{-1mm}
	\caption{Different parabolas}
	\label{Figure1}
\end{figure}

In this setting, the solutions of the quadratic equation $ ax^2+bx+c=0$ are exactly  the $x$-intercept points of the graph of $f(x)$. Depending on the discriminant of the equation, the quadratic function might have two, one or no $x$-intercept points (see \cref{Figure1}).

\subsection{History}
Quadratic equations appear in many Babylonian mathematical texts. It has been suggested that some  date to the Third dynasty of Ur  (circa 2000 BC). The Babylonian scribes were familiar with neither zero nor negative numbers, so the standard form of quadratic equations they usually considered have the following form:
\begin{equation}\label{equ-c}
	ax^2\pm bx=c,
\end{equation}
where $a,b,c$ are non-zero positive numbers.  Based on mathematical texts in the Babylonian mathematical corpus, the algebraic method  Babylonian scribes used to solve the  quadratic equations  $ ax^2\pm bx=c$ was the standard method of completing the square. It should be noted that the Babylonian scribes only considered the positive solution for a quadratic equation. 

In many   Babylonian texts, a quadratic system of equations like \cref{eq-ca} were treated. The Babylonian scribes used the method we explained in \cref{GFS} and  \cref{QSE} to solve such equations. Although in Babylonian mathematics   the  algebraic methods were used to complete the square,  geometric explanations for this method were  provided  later by   the Greek and Muslim mathematicians. Perhaps    Euclid might have been the one who  sparked other mathematicians' curiosity to tackle this issue, as he used a geometric argument  to show that (see \cite{Bur06})\\

\noindent
``\textit{If a straight line is bisected and produced to any point, then the rectangle  contained by the whole (with the added straight line) together with the square  on half the line bisected is equal to the square  on the straight line made up of the half and the part added}''.\\

\begin{figure}[H]
	\centering
	\includegraphics[scale=1]{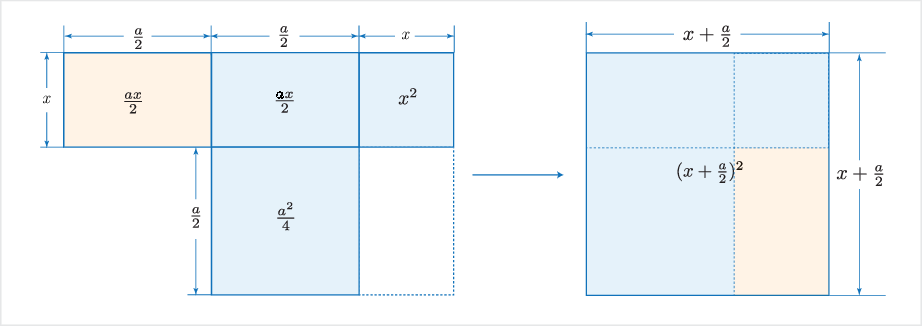}
	\caption{Euclid's completing the square}
	\label{Figure2}
\end{figure}

We have pictured this situation in   \cref{Figure2}. It is evident from       the figure that  the total area  of the colored rectangular regions on the left-hand side    is 
\[x^2+\frac{ax}{2}+\frac{ax}{2}+\left(\frac{a}{2}\right)^2=x^2+ax+ \frac{a^2}{4}.\]
On the other hand, by cutting the orange rectangle  with sides $x$ and $\frac{a}{2}$  and attach it to  the empty rectangle  with the same dimensions, we get a square  of side $x+\frac{a}{2}$ whose area is exactly the same as the sum of areas of colored rectangles. In other words, we have
\[  x^2+ax+ \frac{a^2}{4}=\left(x+\frac{a}{2}\right)^2.\] 
 
Mathematicians after Euclid also tried to use similar  geometric approaches   for solving quadratic equations. One of the greatest Muslim mathematicians   who gave a geometric treatment of this topic  was the prominent Persian mathematician    al-Khwarizmi  (circa 780--850).  He gave two different geometric solutions using completing the square  in order to solve the specific quadratic equation  $ x^2+10x=39$. We explain  his  two methods   in the following.

\begin{figure}[H]
	\centering
	\includegraphics[scale=1]{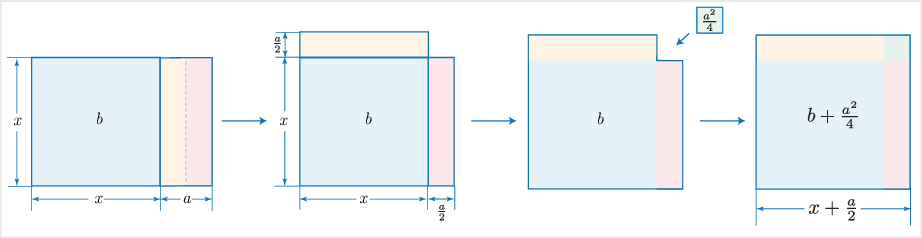}
	\caption{Al-Khwarizmi's first method for completing the square}
	\label{Figure3}
\end{figure}

\noindent
\underline{Al-Khwarizmi's First Method for Completing the Square}\\
The geometric explanation of the  first method  is shown in    \cref{Figure3}. Consider a rectangle with sides $x$ and $x+a$ whose area is the known value $b$ and then slice it off in order to produce  a square with side $x$ and a  rectangle  with sides $x$ and $a$. Next,  cut   off the obtained rectangle in half along the  dotted  red line to get two equal rectangles  with sides $x$ and  $\frac{a}{2}$. Then,  attach these two obtained rectangles along their sides   $x$ to  the top side and   the right  side of the square with side $x$. Finally,  add a small square  with side $\frac{a}{2} $ (and the area $ \frac{a^2}{4}$) to get a complete square  with side $x+\frac{a}{2}$. The area of the resulted square  is exactly the sum of our area   $ b $ and the extra area   $\frac{a^2}{4}$. In other words,  
\[ \left(x+\frac{a}{2}\right)^2-\frac{a^2}{4}=x^2+ax. \]

\begin{figure}[H]
	\centering
	\includegraphics[scale=1]{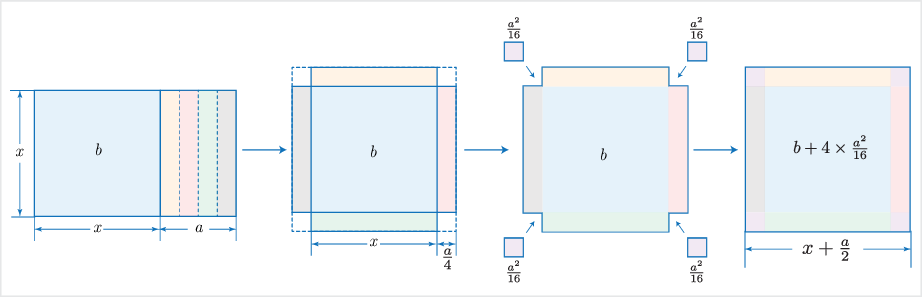}
	\caption{Al-Khwarizmi's second method for  completing the square}
	\label{Figure4}
\end{figure}

	\noindent
\underline{Al-Khwarizmi's Second Method for Completing the Square}\\
In the second method, al-Khwarizmi carries out the procedure as  shown in    \cref{Figure4}. Remove the rectangle  with sides $a$ and $x$ and cut it off along three  dotted red lines in order to get  four equal rectangles   with sides $\frac{a}{4}$ and $x $. 	Then, attach four obtained rectangles   along their sides $x$ to four sides of the square  with side $x$. Finally,   add four small squares   with sides $ \frac{a}{4}$ and areas $ \frac{a^2}{16}$ to get a complete square  with side $x+\frac{a}{2}$. The area of the resulted square is exactly the sum of our area $b$ and the extra four areas $\frac{a^2}{16} $. In other words,
\[ \left(x+\frac{a}{2}\right)^2-4\times\left(\frac{a^2}{16}\right)=x^2+ax. \] 
Note that in both methods we have geometrically verified the algebraic identity
\[ \left(x+\frac{a}{2}\right)^2=x^2+ax+\frac{a^2}{4}. \]

As  Katz mentions in  \cite{Kat09},  Al-Khwarizmi's verbal descriptions of his  methods   are essentially the same as those of the Elamite and Babylonian scribes.   Geometric descriptions of his methods might  have been taken over from his  true  ancestors in the Near East!

\section{Quadratic Equations in SMT}

\subsection{SMT No.\,5-Part 6:   $ ax^2=b$}\label{SS-P6-SMT5}  
In this subsection, simple quadratic equations\index{quadratic equation} are presented without their solutions. The translations of two lines 1 and 2 are 

\noindent
(L1)  50,25 is the area\index{area of a square} (of a square\index{square}). What is the side?\\
(L2)  $ \frac{1}{11}  $ of the area\index{area of a square} (of a square\index{square}) is 4,35. What is the side?

\subsubsection*{Transliteration of Part 6}\label{SSS-P6TI-SMT5}   
\begin{note1} 
	\underline{Obverse III: Lines 1-10}

	(L1)\hspace{2mm} 50,25 a-\v{s}\`{a} n\'{i}gin \textit{mi-nu}

	(L2)\hspace{2mm} \textbf{11} a-\v{s}\`{a} 4,35 n\'{i}gin \textit{mi}-[\textit{nu}]

	(L3)\hspace{2mm} 2 \textbf{11} a-\v{s}\`{a} 9,10 n\'{i}gin \textit{mi-n}[\textit{u}]

	(L4)\hspace{2mm} 1,41,40,25 a-\v{s}\`{a} n\'{i}gin \textit{mi-nu}

	(L5)\hspace{2mm}  \textbf{11 11} a-\v{s}\`{a} 50,25 n\'{i}gin \textit{mi-nu}

	(L6)\hspace{2mm} 2 \textbf{11 11} a-\v{s}\`{a} 1,40,50 n\'{i}gin \textit{mi-}[\textit{nu}]

	(L7)\hspace{2mm} 41,10,25 a-\v{s}\`{a} n\'{i}gin \textit{mi-n}[\textit{u}]

	(L8)\hspace{2mm} \textbf{11 7} a-\v{s}\`{a} 32,5 n\'{i}gin \textit{mi-nu}

	(L9)\hspace{2mm}  2 \textbf{11 7} a-\v{s}\`{a} 1,4,10 n\'{i}gin \textit{mi-nu}

	(L10)  2,44,41,40 a-\v{s}\`{a} n\'{i}gin \textit{mi-nu} 
\end{note1}

\subsubsection*{Mathematical Calculations of Part 6}\label{SSS-P6MC-SMT5}

\noindent
(L1)\hspace{6mm}  $ \displaystyle x^2=50,25$\ \ \ and\ \ \ $x=55$.

\vspace{2mm}
\noindent
(L2)\hspace{6mm}  $ \displaystyle \frac{1}{11}  x^2=4,35$\ \ \ and\ \ \ $x=55$.

\vspace{2mm}
\noindent
(L3)\hspace{6mm}  $ \displaystyle \left(2\times \frac{1}{11}\right) x^2=9,10$\ \ \ and\ \ \ $x=55$.

\vspace{2mm}
\noindent
(L4)\hspace{6mm}  $ \displaystyle x^2=1,41,40,50$\ \ \ and\ \ \ $x=10,5$.

\vspace{2mm}
\noindent
(L5)\hspace{6mm}  $ \displaystyle \left(\frac{1}{11}\times \frac{1}{11}\right) x^2=50,25$\ \ \ and\ \ \ $x=10,5$.

\vspace{2mm}
\noindent
(L6)\hspace{6mm}  $ \displaystyle \left(2\times \frac{1}{11}\times \frac{1}{11}\right) x^2=1,40,50$\ \ \ and\ \ \ $x=10,5$.

\vspace{2mm}
\noindent
(L7)\hspace{6mm}  $ \displaystyle x^2=41,10,25$\ \ \ and\ \ \ $x=6,25$.

\vspace{2mm}
\noindent
(L8)\hspace{6mm}  $ \displaystyle \left(\frac{1}{11}\times \frac{1}{7}\right) x^2=32,5$\ \ \ and\ \ \ $x=6,25$.

\vspace{2mm}
\noindent
(L9)\hspace{6mm}  $ \displaystyle \left(2\times \frac{1}{11}\times \frac{1}{7}\right) x^2=1,4,10$\ \ \ and\ \ \ $x=6,25$.

\vspace{2mm}
\noindent
(L10)\hspace{5mm}  $ \displaystyle x^2=2,44,41,40$\ \ \ and\ \ \ $x=12,50$.

\subsection{SMT No.\,5-Part 7:   $ (ax)^2=b$}\label{SS-P7-SMT5}   
In this subsection too, simple quadratic equations\index{quadratic equation} are presented without their solutions, but the wording is ambiguous as the literal translation of line 19 shows 

\noindent
(L19) 35, the side of a square\index{square}, the area\index{area of a square}, $ \frac{1}{7} $, what is the length\index{length}?

\noindent
With all data available, we cannot help judging that it means 
\[``35^2 = \left(\frac{1}{7} x\right)^2,\  \mathrm{what\  is}\  x ?"\]

\subsubsection*{Transliteration of Part 7}\label{SSS-P7TI-SMT5}   

\begin{note1} 
	\underline{Obverse III: Lines 11-30}

	(L11)\hspace{2mm} 31(sic) n\'{i}gin a-\v{s}\`{a} 2 u\v{s} n\'{i}gin(sic) \textit{mi}-[\textit{nu}]

	(L12)\hspace{2mm} 31(sic) n\'{i}gin a-\v{s}\`{a} 3 u\v{s} [\textit{mi-nu}]

	(L13)\hspace{2mm} 30 n\'{i}gin a-\v{s}\`{a} 4 [u\v{s}  \textit{mi-nu}]

	(L14)\hspace{2mm} 30 n\'{i}gin a-\v{s}\`{a} [$ \frac{\text{2}}{\text{3}} $ u\v{s} \textit{mi-nu}]

	(L15)\hspace{2mm} 30 n\'{i}gin a-\v{s}\`{a} [$ \frac{\text{1}}{\text{2}} $ u\v{s} \textit{mi-nu}]

	(L16)\hspace{2mm} 30 n\'{i}gin [a-\v{s}\`{a} $\frac{\text{1}}{\text{3}} $ u\v{s} \textit{mi-nu}]

	(L17)\hspace{2mm} [30] n\'{i}gin a-\v{s}\`{a} \textbf{4} u\v{s} [\textit{mi-nu}]\\
	(L18)\hspace{2mm} [3]1(sic) n\'{i}gin a-\v{s}\`{a} $\frac{\text{1}}{\text{3}}$ \textbf{4} u\v{s} \textit{mi}-[\textit{nu}]\\
	(L19)\hspace{2mm} 35 n\'{i}gin a-\v{s}\`{a} \textbf{7} u\v{s} \textit{mi-nu}\\
	(L20)\hspace{2mm} 35 n\'{i}gin a-\v{s}\`{a} 2 \textbf{7} u\v{s} \textit{mi-nu}\\
	(L21)\hspace{2mm} 4,5 n\'{i}gin a-\v{s}\`{a} \textbf{7 7} u\v{s} \textit{mi-nu}\\
	(L22)\hspace{2mm} 4,5 n\'{i}gin a-\v{s}\`{a} 2 \textbf{7 7} u\v{s} \textit{mi}-[\textit{nu}]\\
	(L23)\hspace{2mm} 55 n\'{i}gin a-\v{s}\`{a} \textbf{11} u\v{s} \textit{mi}-[\textit{nu}]\\
	(L24)\hspace{2mm} 55 n\'{i}gin a-\v{s}\`{a} 2 \textbf{11} u\v{s} \textit{mi-nu}\\
	(L25)\hspace{2mm} 10,5 n\'{i}gin a-\v{s}\`{a} \textbf{11 11} u\v{s} \textit{mi-nu}\\
	(L26)\hspace{2mm} 10,5 n\'{i}gin a-\v{s}\`{a} 2 \textbf{11 11} u\v{s} \textit{mi-nu}\\
	(L27)\hspace{2mm} 6,25 n\'{i}gin a-\v{s}\`{a} \textbf{11 7} u\v{s} \textit{mi-nu}\\
	(L28)\hspace{2mm} 6,25 n\'{i}gin a-\v{s}\`{a} 2 \textbf{11 8}(sic) u\v{s} \textit{mi-nu}\\
	(L29)\hspace{2mm} 12,50 n\'{i}gin a-\v{s}\`{a} $ \frac{\text{2}}{\text{3}}~\frac{\text{1}}{\text{2}}~\frac{\text{1}}{\text{3}} $   \textbf{11 7} u\v{s} \textit{mi-nu}\\
	(L30)\hspace{2mm} 12,50 n\'{i}gin 2(sic) a-\v{s}\`{a} $ \frac{\text{2}}{\text{3}}~\frac{\text{1}}{\text{2}}~\frac{\text{1}}{\text{3}}$  \textbf{11 7} u\v{s} \textit{mi-nu}
\end{note1}

\subsubsection*{Mathematical Calculations of Part 7}\label{SSS-P7MC-SMT5}

\noindent
(L11)\hspace{6mm}  $ \displaystyle (2x)^2=30^2 $\ \ \ and\ \ \ $x=15$.

\vspace{2mm}
\noindent
(L12)\hspace{6mm}  $ \displaystyle (3x)^2=30^2 $\ \ \ and\ \ \ $x=10$.

\vspace{2mm}
\noindent
(L13)\hspace{6mm}  $ \displaystyle (4x)^2=30^2 $\ \ \ and\ \ \ $x=7;30$.

\vspace{2mm}
\noindent
(L14)\hspace{6mm}  $ \displaystyle \left(\frac{2}{3}x\right)^2=30^2 $\ \ \ and\ \ \ $x=45$.

\vspace{2mm}
\noindent
(L15)\hspace{6mm}  $ \displaystyle \left(\frac{1}{2}x\right)^2=30^2 $\ \ \ and\ \ \ $x=1,0$.

\vspace{2mm}
\noindent
(L16)\hspace{6mm}  $ \displaystyle \left(\frac{1}{3}x\right)^2=30^2 $\ \ \ and\ \ \ $x=1,30$.

\vspace{2mm}
\noindent
(L17)\hspace{6mm}  $ \displaystyle \left(\frac{1}{4}x\right)^2=30^2 $\ \ \ and\ \ \ $x=2,0$.

\vspace{2mm}
\noindent
(L18)\hspace{6mm}  $ \displaystyle \left[ \left(\frac{1}{3}\times \frac{1}{4}\right) x\right]^2=30^2 $\ \ \ and\ \ \ $x=6,0$.

\vspace{2mm}
\noindent
(L19)\hspace{6mm}  $ \displaystyle \left(\frac{1}{7}x\right)^2=35^2 $\ \ \ and\ \ \ $x=6,5$.

\vspace{2mm}
\noindent
(L20)\hspace{6mm}  $ \displaystyle \left[ \left(2\times\frac{1}{7}\right) x\right]^2=35^2 $\ \ \ and\ \ \ $x=3,2;30$.

\vspace{2mm}
\noindent
(L21)\hspace{6mm}  $ \displaystyle \left[ \left(\frac{1}{7}\times\frac{1}{7}\right)x\right]^2=(4,5)^2 $\ \ \ and\ \ \ $x=3,20,5$.

\vspace{2mm}
\noindent
(L22)\hspace{6mm}  $ \displaystyle \left[ \left(2\times\frac{1}{7}\times\frac{1}{7}\right) x\right]^2=(4,5)^2 $\ \ \ and\ \ \ $x=1,40,2;30$.

\vspace{2mm}
\noindent
(L23)\hspace{6mm}  $ \displaystyle \left( \frac{1}{11}  x\right)^2=55^2 $\ \ \ and\ \ \ $x=10,5$.

\vspace{2mm}
\noindent
(L24)\hspace{6mm}  $ \displaystyle \left[ \left(2\times\frac{1}{11} \right) x\right]^2=55^2 $\ \ \ and\ \ \ $x=5,2;30$.

\vspace{2mm}
\noindent
(L25)\hspace{6mm}  $ \displaystyle \left[ \left( \frac{1}{11}\times\frac{1}{11}\right) x\right]^2=(10,5)^2 $\ \ \ and\ \ \ $x=20,20,5$.

\vspace{2mm}
\noindent
(L26)\hspace{6mm}  $ \displaystyle \left[ \left(2\times\frac{1}{11}\times\frac{1}{11}\right) x\right]^2=(10,5)^2 $\ \ \ and\ \ \ $x=10,10,2;30$.

\vspace{2mm}
\noindent
(L27)\hspace{6mm}  $ \displaystyle \left[ \left(\frac{1}{11}\times\frac{1}{7}\right) x\right]^2=(6,25)^2 $\ \ \ and\ \ \ $x=8,14,5$.

\vspace{2mm}
\noindent
(L28)\hspace{6mm}  $ \displaystyle \left[ \left(2\times\frac{1}{11}\times\frac{1}{7}\right) x\right]^2=(6,25)^2 $\ \ \ and\ \ \ $x=4,7,2;30$.

\vspace{2mm}
\noindent
(L29)\hspace{6mm}  $ \displaystyle \left[ \left(\frac{2}{3} \times \frac{1}{2} \times \frac{1}{3}\times\frac{1}{11}\times\frac{1}{7}\right) x\right]^2=(12,50)^2 $\ \ \ and\ \ \ $x=2,28,13,30$.

\vspace{2mm}
\noindent
(L30)\hspace{6mm}  $ \displaystyle \left[ \left(2\times \frac{2}{3} \times \frac{1}{2} \times \frac{1}{3}\times\frac{1}{11}\times\frac{1}{7}\right) x\right]^2=(12,50)^2 $\ \ \ and\ \ \ $x=1,14,6,45$.

\subsection{SMT No.\,5-Part 8:   $ x^2+(ax)^2 =b$}\label{SS-P8-SMT5} 
In this subsection, simple quadratic equations\index{quadratic equation} are    presented. The translation of line 31 is

\noindent
(L31) The square (literally: area\index{area of a square}) of twice the length\index{length} and the square (of the length\index{length}) are added up to 1,15,0. (What is the length\index{length}?)

\subsubsection*{Transliteration of Part 8}\label{SSS-P8TI-SMT5}   
\begin{note1} 
	\underline{Obverse III: Lines 31-46} \\
	(L31)\hspace{2mm} a-\v{s}\`{a} 2 u\v{s} \textit{\`{u}} a-\v{s}\`{a} gar-\textit{ma} 1,15\\
	(L32)\hspace{2mm} a-\v{s}\`{a} 3 u\v{s} \textit{\`{u}} a-\v{s}\`{a} gar-gar-\textit{ma} 2,[30]\\
	(L33)\hspace{2mm} a-\v{s}\`{a} 4 u\v{s} \textit{\`{u}} a-\v{s}\`{a} gar-gar-\textit{ma} 4,1[5]\\
	(L34)\hspace{2mm} a-\v{s}\`{a} \textit{\`{u}} a-\v{s}\`{a} $ \frac{\text{2}}{\text{3}} $ u\v{s} gar-gar-\textit{ma} 21,40\\
	(L35)\hspace{2mm} a-\v{s}\`{a} \textit{\`{u}} a-\v{s}\`{a} $ \frac{\text{1}}{\text{2}} $ u\v{s} gar-gar-\textit{ma} 1[8],45\\
	(L36)\hspace{2mm} a-\v{s}\`{a} \textit{\`{u}} a-\v{s}\`{a} $\frac{\text{1}}{\text{3}} $ u\v{s} gar-gar-\textit{ma} 16,40\\
	(L37)\hspace{2mm} a-\v{s}\`{a} \textit{\`{u}} a-\v{s}\`{a} \textbf{4} u\v{s} gar-gar-\textit{ma} 15,56,15\\
	(L38)\hspace{2mm} a-\v{s}\`{a} \textit{\`{u}} a-\v{s}\`{a} $ \frac{\text{1}}{\text{3}} $ \textbf{4} u\v{s} gar-gar-\textit{ma} 15,6,15\\ 
	(L39)\hspace{2mm} a-\v{s}\`{a} \textit{\`{u}} a-\v{s}\`{a} \textbf{7} u\v{s} gar-gar-\textit{ma} 20,50\\
	(L40)\hspace{2mm} a-\v{s}\`{a} \textit{\`{u}} a-\v{s}\`{a} 2 \textbf{7} u\v{s} gar-gar-\textit{ma} 22,5\\
	(L41)\hspace{2mm} a-\v{s}\`{a} \textit{\`{u}} a-\v{s}\`{a} \textbf{7 7} u\v{s} gar-gar-\textit{ma} 16,40,50\\
	(L42)\hspace{2mm} a-[\v{s}\`{a} \textit{\`{u}} a-\v{s}\`{a}] 2 \textbf{7 7} u\v{s} gar-gar-\textit{ma} 16,42,5\\
	(L43)\hspace{2mm} [a-\v{s}\`{a}] \textit{\`{u}} a-\v{s}\`{a} \textbf{11} u\v{s} gar-gar-\textit{ma} [50],50\\
	(L44)\hspace{2mm} [a-\v{s}\`{a} \textit{\`{u}} a-\v{s}\`{a}] 2 \textbf{11} u\v{s} gar-gar-\textit{ma} 52,5\\
	(L45)\hspace{2mm} [a-\v{s}\`{a} \textit{\`{u}} a-\v{s}\`{a} \textbf{11 11} u\v{s} gar-gar-\textit{ma} 1,4]1,40,50\\
	(L46)\hspace{2mm} [a-\v{s}\`{a} \textit{\`{u}} a-\v{s}\`{a} 2 \textbf{11 11} u\v{s} gar-gar-\textit{ma} 1,41,42,5]
\end{note1}

\subsubsection*{Mathematical Calculations of Part 8}\label{SSS-P8MC-SMT5}

\noindent
(L31)\hspace{6mm}  $ \displaystyle  x^2+(2x)^2 =1,15,0 $\ \ \ and\ \ \ $x=30$.

\vspace{2mm}
\noindent
(L32)\hspace{6mm}  $ \displaystyle  x^2+(3x)^2 =2,30,0 $\ \ \ and\ \ \ $x=30$.

\vspace{2mm}
\noindent
(L33)\hspace{6mm}  $ \displaystyle x^2+(4x)^2 =4,15,0 $\ \ \ and\ \ \ $x=30$.

\vspace{2mm}
\noindent
(L34)\hspace{6mm}  $ \displaystyle x^2+\left(\frac{2}{3}x\right)^2=12,40 $\ \ \ and\ \ \ $x=30$.

\vspace{2mm}
\noindent
(L35)\hspace{6mm}  $ \displaystyle x^2+\left(\frac{1}{2}x\right)^2=18,45 $\ \ \ and\ \ \ $x=30$.

\vspace{2mm}
\noindent
(L36)\hspace{6mm}  $ \displaystyle x^2+\left(\frac{1}{3}x\right)^2=16,40 $\ \ \ and\ \ \ $x=30$.

\vspace{2mm}
\noindent
(L37)\hspace{6mm}  $ \displaystyle x^2+\left(\frac{1}{4}x\right)^2=15,56;15 $\ \ \ and\ \ \ $x=30$.

\vspace{2mm}
\noindent
(L38)\hspace{6mm}  $ \displaystyle x^2+\left[ \left(\frac{1}{3}\times\frac{1}{4}\right)x\right]^2=15,6;15 $\ \ \ and\ \ \ $x=30$.

\vspace{2mm}
\noindent
(L39)\hspace{6mm}  $ \displaystyle x^2+\left(\frac{1}{7}x\right)^2=20,50 $\ \ \ and\ \ \ $x=35$.

\vspace{2mm}
\noindent
(L40)\hspace{6mm}  $ \displaystyle x^2+\left[ \left(2\times \frac{1}{7}\right)x\right]^2=22,5 $\ \ \ and\ \ \ $x=35$.

\vspace{2mm}
\noindent
(L41)\hspace{6mm}  $ \displaystyle x^2+\left[ \left(\frac{1}{7}\times \frac{1}{7}\right)x\right]^2=16,40,5 $\ \ \ and\ \ \ $x=4,5$.

\vspace{2mm}
\noindent
(L42)\hspace{6mm}  $ \displaystyle x^2+\left[ \left(2\times \frac{1}{7}\times \frac{1}{7}\right)x\right]^2=16,42,5 $\ \ \ and\ \ \ $x=4,5$.

\vspace{2mm}
\noindent
(L43)\hspace{6mm}  $ \displaystyle x^2+\left( \frac{1}{11} x\right)^2=50,50 $\ \ \ and\ \ \ $x=55$.

\vspace{2mm}
\noindent
(L44)\hspace{6mm}  $ \displaystyle x^2+\left[ \left(2\times \frac{1}{11}\right)x\right]^2=52,5 $\ \ \ and\ \ \ $x=55$.

\vspace{2mm}
\noindent
(L45)\hspace{6mm}  $ \displaystyle x^2+\left[ \left(\frac{1}{11}\times \frac{1}{11}\right)x\right]^2=1,41,40,50 $\ \ \ and\ \ \ $x=10,5$.

\vspace{2mm}
\noindent
(L46)\hspace{6mm}  $ \displaystyle x^2+\left[ \left(2\times\frac{1}{11}\times \frac{1}{11}\right)x\right]^2=1,41,42,5 $\ \ \ and\ \ \ $x=10,5$.

\subsection{SMT No.\,5-Part 9:   $ x^2-(ax)^2 =b$}\label{SS-P9-SMT5} 
This subsection deals with quadratic equations\index{quadratic equation} which are  conjugate\index{conjugate quadratic equations} to those in the previous subsection (see   \cref{SS-P8-SMT5}). The statement of the problem is as follows:

\noindent
(L7) The square (of a length\index{length}) exceeds the square of $ \frac{1}{7} $ of the length\index{length} by 20,0. (What is the length\index{length}?)

\subsubsection*{Transliteration of Part 9}\label{SSS-P9TI-SMT5}   
\begin{note1} 
	\underline{Reverse I: Lines 2-20}\\
	(L2)\hspace{3mm} [a-\v{s}\`{a} ugu a-\v{s}\`{a}] $ \frac{\text{2}}{\text{3}} $ u\v{s} 8,20 dir[ig]\\
	(L3)\hspace{3mm} [a-\v{s}\`{a} ugu a-\v{s}\`{a} $ \frac{\text{1}}{\text{2}} $] u\v{s} 11,15 dirig\\
	(L4)\hspace{3mm} [a-\v{s}\`{a} ugu a-\v{s}\`{a} $ \frac{\text{1}}{\text{3}} $] u\v{s} 13,20 $<$dirig$>$\\
	(L5)\hspace{3mm} [a-\v{s}\`{a} ugu a-\v{s}\`{a} \textbf{4}] u\v{s} 14,3,45 $<$dirig$>$\\
	(L6)\hspace{3mm} a-\v{s}\`{a} ugu a-\v{s}\`{a} [$ \frac{\text{1}}{\text{3}} $ \textbf{4}] u\v{s} 14,53,45 $<$dirig$>$\\
	(L7)\hspace{3mm} a-\v{s}\`{a} ugu a-\v{s}\`{a} \textbf{7} u\v{s} 20 dirig\\
	(L8)\hspace{3mm} a-\v{s}\`{a} ugu a-\v{s}\`{a} 2 \textbf{7} u\v{s} 18,45 $<$dirig$>$\\
	(L9)\hspace{3mm} a-\v{s}[\`{a}] ug[u] a-\v{s}\`{a} \textbf{7 7} u\v{s} 16,40 dirig\\
	(L10)\hspace{1mm} a-\v{s}\`{a} ugu a-\v{s}\`{a} 2 \textbf{7 7} u\v{s} 16,42(sic),5(sic) $<$dirig$>$\\
	(L11)\hspace{1mm} [a-\v{s}\`{a}] ugu a-\v{s}\`{a} \textbf{11} u\v{s} 50 [dirig]\\
	(L12)\hspace{1mm} [a-\v{s}\`{a}] ugu a-\v{s}\`{a} 2 \textbf{11} u\v{s} 48,4[5] [dir]ig\\
	(L13)\hspace{1mm} a-\v{s}\`{a} ugu a-\v{s}\`{a} \textbf{11 11} u\v{s} 1,41,40 [dirig]\\
	(L14)\hspace{1mm} a-\v{s}\`{a} ugu a-\v{s}\`{a} 2 \textbf{11 11} u\v{s} 1,41,38,[4]5 $<$dirig$>$\\
	(L15)\hspace{1mm} a-\v{s}\`{a} ugu a-\v{s}\`{a} \textbf{11 7} u\v{s} 41,10 dirig\\
	(L16)\hspace{1mm} a-\v{s}\`{a} ugu a-\v{s}\`{a} 2 \textbf{11 7} u\v{s} $<$4$>$1,8,45 dirig\\
	(L17)\hspace{1mm} a-\v{s}\`{a} ugu a-\v{s}\`{a} $ \frac{\text{2}}{\text{3}}~\frac{\text{1}}{\text{2}}~\frac{\text{1}}{\text{3}}$  \textbf{11 7} u\v{s} 2,44,41,30\\
	(L18)\hspace{1mm}  8,45,55,33,20 dirig\\
	(L19)\hspace{1mm} a-\v{s}\`{a} ugu a-\v{s}\`{a} 2 $ \frac{\text{2}}{\text{3}}~\frac{\text{1}}{\text{2}}~\frac{\text{1}}{\text{3}} $ \textbf{11 7} u\v{s} 2,44,41,30\\
	(L20)\hspace{1mm}  5,3,42,13,20 dirig
\end{note1}

\subsubsection*{Mathematical Calculations of Part 9}\label{SSS-P9MC-SMT5}

\noindent
(L2)\hspace{6mm}  $ \displaystyle x^2-\left(\frac{2}{3}x \right)^2 =8,20 $\ \ \ and\ \ \ $x=30$.

\vspace{2mm}
\noindent
(L3)\hspace{6mm}  $ \displaystyle x^2-\left(\frac{1}{2}x \right)^2 =11,15 $\ \ \ and\ \ \ $x=30$.  

\vspace{2mm}
\noindent
(L4)\hspace{6mm}  $ \displaystyle x^2-\left(\frac{1}{3}x \right)^2 =13,20 $\ \ \ and\ \ \ $x=30$.  

\vspace{2mm}
\noindent
(L5)\hspace{6mm}  $ \displaystyle x^2-\left(\frac{1}{4}x \right)^2 =14,3;45 $\ \ \ and\ \ \ $x=30$.  

\vspace{2mm}
\noindent
(L6)\hspace{6mm}  $ \displaystyle x^2-\left[ \left(\frac{1}{3}\times \frac{1}{4}\right)x \right]^2 =14,53;45 $\ \ \ and\ \ \ $x=30$.  

\vspace{2mm}
\noindent
(L7)\hspace{6mm}  $ \displaystyle x^2-\left(\frac{1}{7}x \right)^2 =20,0 $\ \ \ and\ \ \ $x=35$.  

\vspace{2mm}
\noindent
(L8)\hspace{6mm}  $ \displaystyle x^2-\left[ \left(2\times \frac{1}{7}\right)x \right]^2 =18,45 $\ \ \ and\ \ \ $x=35$.

\vspace{2mm}
\noindent
(L9)\hspace{6mm}  $ \displaystyle x^2-\left[ \left(\frac{1}{7}\times \frac{1}{7}\right)x \right]^2 =16,40,0$\ \ \ and\ \ \ $x=4,5$.  

\vspace{2mm}
\noindent
(L10)\hspace{6mm}  $ \displaystyle x^2-\left[ \left(2\times\frac{1}{7}\times \frac{1}{7}\right)x \right]^2 =16,38,45 $\ \ \ and\ \ \ $x=4,5$.  

\vspace{2mm}
\noindent
(L11)\hspace{6mm}  $ \displaystyle x^2-\left(  \frac{1}{11}\ x \right)^2 =50,0 $\ \ \ and\ \ \ $x=55$.

\vspace{2mm}
\noindent
(L12)\hspace{6mm}  $ \displaystyle x^2-\left[ \left(2\times \frac{1}{11}\right)x \right]^2 =48,45 $\ \ \ and\ \ \ $x=55$.  

\vspace{2mm}
\noindent
(L13)\hspace{6mm}  $ \displaystyle x^2-\left[ \left(\frac{1}{11}\times \frac{1}{11}\right)x \right]^2 =1,41,40,0$\ \ \ and\ \ \ $x=10,5$.  

\vspace{2mm}
\noindent
(L14)\hspace{6mm}  $ \displaystyle x^2-\left[ \left(2\times\frac{1}{11}\times \frac{1}{11}\right)x \right]^2 =1,41,38,45 $\ \ \ and\ \ \ $x=10,5$.  

\vspace{2mm}
\noindent
(L15)\hspace{6mm}  $ \displaystyle x^2-\left[ \left(\frac{1}{11}\times \frac{1}{7}\right)x \right]^2 =41,10,0 $\ \ \ and\ \ \ $x=6,25$.  

\vspace{2mm}
\noindent
(L16)\hspace{6mm}  $ \displaystyle x^2-\left[ \left(2\times\frac{1}{11}\times \frac{1}{7}\right)x \right]^2 =41,8,45 $\ \ \ and\ \ \ $x=6,25$.

\vspace{2mm}
\noindent
(L17-18)\hspace{0mm}  $ \displaystyle x^2-\left[ \left(\frac{2}{3}\times\frac{1}{2}\times \frac{1}{3}\times\frac{1}{11}\times \frac{1}{7}\right)x \right]^2 =2,44,41,38;45,55,33,20 $\ \ \ and\linebreak
\\
${}$ \hspace{16mm} $x=12,50$.

\vspace{2mm}
\noindent
(L19-20)\hspace{0mm}  $ \displaystyle x^2-\left[ \left(2\times\frac{2}{3}\times\frac{1}{2}\times \frac{1}{3}\times\frac{1}{11}\times \frac{1}{7}\right)x \right]^2 =2,44,41,35;3,42,13,20 $\ \ \ and   \linebreak
\\
$ {}$\hspace{17mm} $x=12,50$.

\subsection{SMT No.\,5-Part 10:    $ x^2+ax =b$}\label{SS-P10-SMT5} 
In this subsection, a new group of  quadratic equations  are    presented which can be solved    by completing the square. 
The translation of the first problem in line 21  is as follows: 

\noindent
(L21) To the square of my side one length  (my side) is added, and (the result is) 0;45. (What is my side?).

This means $x^2 + x = 0;45$  in    modern mathematical notation. In the following problems, the term ``To the square of my side'' is omitted. 

It is   interesting   that the same quadratic equation  as $x^2 + x = 0;45$  occurs in the first problem of the mathematical tablet  \textbf{BM 13901}\footnote{This tablet  is one of the oldest known mathematical texts held in the British Museum. It contains about twenty-four problems and their solutions.}    which also deals with quadratic equations  and  uses completing the square  (see \cite{Mur03-1}, for more information).

\subsubsection*{Transliteration of Part 10}\label{SSS-P10TI-SMT5}   

\begin{note1} 
	\underline{Reverse I: Lines 21-47} \\
	(L21)\hspace{3mm} \textit{a-na} a-\v{s}\`{a} nigin-\textit{ia} 1 u\v{s} dah-\textit{ma} 45\\
	(L22)\hspace{3mm} 2 u\v{s} [da]h-\textit{ma} 1,15\\
	(L23)\hspace{3mm} 3 u\v{s} [da]h-\textit{ma} 1,45\\
	(L24)\hspace{3mm} \textbf{4} u\v{s} dah-\textit{ma} 2,15\\
	(L25)\hspace{3mm} $ \frac{\text{2}}{\text{3}} $ u\v{s} dah-\textit{ma} 35\\
	(L26)\hspace{3mm} $ \frac{\text{1}}{\text{2}} $ u\v{s} dah-\textit{ma} 30\\
	(L27)\hspace{3mm} $ \frac{\text{1}}{\text{3}} $ u\v{s} dah-\textit{ma} 25\\
	(L28)\hspace{3mm} \textbf{4} u\v{s} dah-\textit{ma} 22,30\\
	(L29)\hspace{3mm} $ \frac{\text{1}}{\text{3}} $ \textbf{4} u\v{s} dah-\textit{ma} 17,30\\
	(L30)\hspace{3mm} \textbf{7} u\v{s} dah-\textit{ma} 25,2[5]\\
	(L31)\hspace{3mm} 2 \textbf{7} u\v{s} dah-\textit{ma} 30,2[5]\\
	(L32)\hspace{3mm} \textbf{7 7} u\v{s} dah-\textit{ma} [16,45,25]\\
	(L33)\hspace{3mm} 2 \textbf{7 7} u\v{s} dah-\textit{ma} [16,50,25]\\
	(L34)\hspace{3mm} 1 $ \frac{\text{2}}{\text{3}} $ u\v{s} dah-\textit{ma} [1,5]\\
	(L35)\hspace{3mm} 1 $ \frac{\text{1}}{\text{2}}$ u\v{s} dah-\textit{ma} 1\\
	(L36)\hspace{3mm} 1 $ \frac{\text{1}}{\text{3}} $ u\v{s} dah-\textit{ma} 55\\
	(L37)\hspace{3mm} 1 \textbf{4} u\v{s} dah-\textit{ma} 52,30\\
	(L38)\hspace{3mm} 1 $ \frac{\text{1}}{\text{3}} $ \textbf{4} u\v{s} dah-\textit{ma} 47,30\\
	(L39)\hspace{3mm} 1 \textbf{7} u\v{s} dah-\textit{ma} 1,0,25\\
	(L40)\hspace{3mm} 1 2 \textbf{7} u\v{s} dah-\textit{ma} 1,5,25 \\
	(L41)\hspace{3mm} 1 \textbf{7 7} u\v{s} dah-\textit{ma} 20,45(sic),25\\
	(L42)\hspace{3mm} 1 2 \textbf{7 7} u\v{s} dah-\textit{ma} 20,50(sic),25\\
	(L43)\hspace{3mm} 2 $ \frac{\text{1}}{\text{2}} $ u\v{s} dah-\textit{ma} 1,30\\
	(L44)\hspace{3mm} 3 $ \frac{\text{1}}{\text{3}} $ u\v{s} dah-\textit{ma} 1,55\\
	(L45)\hspace{3mm} 4 \textbf{4} u\v{s} dah-\textit{ma} 2,22,30\\
	(L46)\hspace{3mm} 7 igi-7 u\v{s} dah-\textit{ma} 24,35\\
	(L47)\hspace{3mm} 7 2 igi-7 u\v{s} dah-\textit{ma} 24,40
\end{note1}

\subsubsection*{Mathematical Calculations of Part 10}\label{SSS-P10MC-SMT5}

\noindent
(L21)\hspace{6mm}  $ \displaystyle x^2 + x=0;45 $\ \ \ and\ \ \ $x=0;30$.

\vspace{2mm}
\noindent
(L22)\hspace{6mm}  $ \displaystyle x^2 + 2x=1;15 $\ \ \ and\ \ \ $x=0;30$.

\vspace{2mm}
\noindent
(L23)\hspace{6mm}  $ \displaystyle x^2 + 3x=1;45 $\ \ \ and\ \ \ $x=0;30$.

\vspace{2mm}
\noindent
(L24)\hspace{6mm}  $ \displaystyle x^2 + 4x=2;15 $\ \ \ and\ \ \ $x=0;30$.

\vspace{2mm}
\noindent
(L25)\hspace{6mm}  $ \displaystyle x^2 + \frac{2}{3}  x=0;35 $\ \ \ and\ \ \ $x=0;30$.

\vspace{2mm}
\noindent
(L26)\hspace{6mm}  $ \displaystyle x^2 + \frac{1}{2}  x=0;30 $\ \ \ and\ \ \ $x=0;30$.

\vspace{2mm}
\noindent
(L27)\hspace{6mm}  $ \displaystyle x^2 + \frac{1}{3}  x=0;25 $\ \ \ and\ \ \ $x=0;30$.

\vspace{2mm}
\noindent
(L28)\hspace{6mm}  $ \displaystyle x^2 + \frac{1}{4}  x=0;22,30 $\ \ \ and\ \ \ $x=0;30$.

\vspace{2mm}
\noindent
(L29)\hspace{6mm}  $ \displaystyle x^2 + \left( \frac{1}{3}\times \frac{1}{4} \right) x=0;17,30 $\ \ \ and\ \ \ $x=0;30$.

\vspace{2mm}
\noindent
(L30)\hspace{6mm}  $ \displaystyle x^2 + \frac{1}{7}  x=0;25,25 $\ \ \ and\ \ \ $x=0;35$.

\vspace{2mm}
\noindent
(L31)\hspace{6mm}  $ \displaystyle x^2 +\left( 2\times \frac{1}{7}\right)  x=0;30,25 $\ \ \ and\ \ \ $x=0;35$.

\vspace{2mm}
\noindent
(L32)\hspace{6mm}  $ \displaystyle x^2 +\left( \frac{1}{7}\times \frac{1}{7}\right)  x=16;45,25 $\ \ \ and\ \ \ $x=4;5$.

\vspace{2mm}
\noindent
(L33)\hspace{6mm}  $ \displaystyle x^2 +\left( 2\times \frac{1}{7}\times \frac{1}{7}\right)  x=16;50,25 $\ \ \ and\ \ \ $x=4;5$.

\vspace{2mm}
\noindent
(L34)\hspace{6mm}  $ \displaystyle x^2 +\left(1+ \frac{2}{3}\right)  x=1;5 $\ \ \ and\ \ \ $x=0;30$.

\vspace{2mm}
\noindent
(L35)\hspace{6mm}  $ \displaystyle x^2 +\left(1+ \frac{1}{2}\right)  x=1  $\ \ \ and\ \ \ $x=0;30$.

\vspace{2mm}
\noindent
(L36)\hspace{6mm}  $ \displaystyle x^2 +\left(1+ \frac{1}{3}\right)  x=0;55 $\ \ \ and\ \ \ $x=0;30$.

\vspace{2mm}
\noindent
(L37)\hspace{6mm}  $ \displaystyle x^2 +\left(1+ \frac{1}{4}\right)  x=0;52,30 $\ \ \ and\ \ \ $x=0;30$.

\vspace{2mm}
\noindent
(L38)\hspace{6mm}  $ \displaystyle x^2 +\left[1+ \left(\frac{1}{3}\times \frac{1}{4}\right)\right]  x=0;47,30 $\ \ \ and\ \ \ $x=0;30$.

\vspace{2mm}
\noindent
(L39)\hspace{6mm}  $ \displaystyle x^2 +\left(1+ \frac{1}{7}\right)  x=1;0,25 $\ \ \ and\ \ \ $x=0;35$.

\vspace{2mm}
\noindent
(L40)\hspace{6mm}  $ \displaystyle x^2 +\left[1+ \left(2\times\frac{1}{7}\right)\right] x=1;5,25 $\ \ \ and\ \ \ $x=0;35$.

\vspace{2mm}
\noindent
(L41)\hspace{6mm}  $ \displaystyle x^2 +\left[1+ \left(\frac{1}{7}\times\frac{1}{7}\right)\right]  x=20;50,25 $\ \ \ and\ \ \ $x=4;5$.

\vspace{2mm}
\noindent
(L42)\hspace{6mm}  $ \displaystyle x^2 +\left[1+ \left(2\times \frac{1}{7}\times\frac{1}{7}\right)\right]  x=20;55,25 $\ \ \ and\ \ \ $x=4;5$.

\vspace{2mm}
\noindent
(L43)\hspace{6mm}  $ \displaystyle x^2 +\left(2+ \frac{1}{2}\right)  x=1;30 $\ \ \ and\ \ \ $x=0;30$.

\vspace{2mm}
\noindent
(L44)\hspace{6mm}  $ \displaystyle x^2 +\left(3+ \frac{1}{3}\right)  x=1;55 $\ \ \ and\ \ \ $x=0;30$.

\vspace{2mm}
\noindent
(L45)\hspace{6mm}  $ \displaystyle x^2 +\left(4+ \frac{1}{4}\right)  x=2;22,30 $\ \ \ and\ \ \ $x=0;30$.

\vspace{2mm}
\noindent
(L46)\hspace{6mm}  $ \displaystyle x^2 +\left(7+ \frac{1}{7}\right)  x=24,35 $\ \ \ and\ \ \ $x=35$.

\vspace{2mm}
\noindent
(L47)\hspace{6mm}  $ \displaystyle x^2 +\left[7+ \left(2\times \frac{1}{7}\right)\right]  x=24,40 $\ \ \ and\ \ \ $x=35$.

\subsection{SMT No.\,5-Part 11-a:     $ x^2-ax =b$}\label{SS-P11a-SMT5} 
This subsection contains quadratic equations which are conjugate  to   those in  the previous subsection (see   \cref{SS-P10-SMT5}), although the first several lines are   lost.

Note that the scribes have implicitly assumed that $x>a$.  The statement of the first lost problem would be:

\noindent
(L1)  From the square of my side one length  (my side) is subtracted, and (the result is) 14,30. (What is my side?)

This gives the quadratic equation  $x^2 - x = 14,30$  which appears  as the second problem in the mathematical tablet  \textbf{BM 13901}.

\subsubsection*{Transliteration of Part 11-a}\label{SSS-P11aTI-SMT5}   
\begin{note1} 
	\underline{Reverse II: Lines 1-17} \\
	(L1)\hspace{3mm} [2 \textbf{7} u\v{s} zi-\textit{ma} 10,2]5\\
	(L2)\hspace{3mm} [\textbf{7 7} u\v{s} zi-\textit{ma} 1]6,35,25\\
	(L3)\hspace{3mm} [2 \textbf{7 7} u\v{s} zi-\textit{m}]\textit{a} 16,30,25\\
	(L4)\hspace{3mm} [1 $ \frac{\text{2}}{\text{3}} $ u\v{s} z]i-\textit{ma} 14,10\\
	(L5)\hspace{3mm} [1 $ \frac{\text{1}}{\text{2}} $ u\v{s}] zi-\textit{ma} 14,15\\
	(L6)\hspace{3mm} [1 $ \frac{\text{1}}{\text{3}} $] u\v{s} zi-\textit{ma} 14,20\\
	(L7)\hspace{3mm} [1] \textbf{4} u\v{s} zi-\textit{ma} 14,22,30\\
	(L8)\hspace{3mm} 1 $ \frac{\text{1}}{\text{3}} $ \textbf{4} u\v{s} zi-\textit{ma} 14,27,30\\
	(L9)\hspace{3mm} 1 \textbf{7} u\v{s} zi-\textit{ma} 19,45\\
	(L10)\hspace{1mm} 1 2 \textbf{7} u\v{s} zi-\textit{ma} 19,40\\
	(L11)\hspace{1mm} 1 \textbf{7 7} u\v{s} zi-\textit{ma} 12,35(sic),25\\
	(L12)\hspace{1mm} [1 2 \textbf{7}] \textbf{7} u\v{s} zi-\textit{ma} 12,30(sic),25\\
	(L13)\hspace{1mm} [2 $ \frac{\text{1}}{\text{2}} $ u\v{s}] zi-\textit{ma} 13,25(sic)\\
	(L14)\hspace{1mm} [3 $ \frac{\text{1}}{\text{3}}$ u\v{s}] zi-\textit{ma} 13,20\\
	(L15)\hspace{1mm} [4 \textbf{4} u]\v{s} zi-\textit{ma} 12,52,30\\
	(L16)\hspace{1mm} [7 igi-7] u\v{s} zi-\textit{ma} 16,15\\
	(L17)\hspace{1mm} [7 2 igi-7] u\v{s} zi-\textit{ma} 16,10
\end{note1}

\subsubsection*{Mathematical Calculations of Part 11-a}\label{SSS-P11aMC-SMT5}

\noindent
(L1)\hspace{6mm}  $ \displaystyle x^2 -\left( 2\times \frac{1}{7} \right) x=0;10,25 $\ \ \ and\ \ \ $x=0;35$.

\vspace{2mm}
\noindent
(L2)\hspace{6mm}  $ \displaystyle x^2 -\left(\frac{1}{7} \times \frac{1}{7} \right) x=16;35,25 $\ \ \ and\ \ \ $x=4;5$.

\vspace{2mm}
\noindent
(L3)\hspace{6mm}  $ \displaystyle x^2 -\left(2\times \frac{1}{7} \times \frac{1}{7} \right) x=16;30,25 $\ \ \ and\ \ \ $x=4;5$.

\vspace{2mm}
\noindent
(L4)\hspace{6mm}  $ \displaystyle x^2 -\left( 1+ \frac{2}{3} \right) x=14,10 $\ \ \ and\ \ \ $x=30$.

\vspace{2mm}
\noindent
(L5)\hspace{6mm}  $ \displaystyle x^2 -\left( 1+ \frac{1}{2} \right) x=14,15 $\ \ \ and\ \ \ $x=30$.

\vspace{2mm}
\noindent
(L6)\hspace{6mm}  $ \displaystyle x^2 -\left( 1+ \frac{1}{3} \right) x=14,20 $\ \ \ and\ \ \ $x=30$.

\vspace{2mm}
\noindent
(L7)\hspace{6mm}  $ \displaystyle x^2 -\left( 1+ \frac{1}{4} \right) x=14,22;30 $\ \ \ and\ \ \ $x=30$.

\vspace{2mm}
\noindent
(L8)\hspace{6mm}  $ \displaystyle x^2 -\left( 1+ \frac{1}{3}+\frac{1}{4} \right) x=14,27;30 $\ \ \ and\ \ \ $x=30$.

\vspace{2mm}
\noindent
(L9)\hspace{6mm}  $ \displaystyle x^2 -\left( 1+  \frac{1}{7} \right) x=19,45 $\ \ \ and\ \ \ $x=35$. 

\vspace{2mm}
\noindent
(L10)\hspace{4mm}  $ \displaystyle x^2 -\left[ 1+ \left(2\times \frac{1}{7}\right) \right] x=19,40 $\ \ \ and\ \ \ $x=35$.

\vspace{2mm}
\noindent
(L11)\hspace{4mm}  $ \displaystyle x^2 -\left[ 1+ \left(\frac{1}{7}\times\frac{1}{7}\right) \right] x=12;30,25 $\ \ \ and\ \ \ $x=4;5$.

\vspace{2mm}
\noindent
(L12)\hspace{4mm}  $ \displaystyle x^2 -\left[ 1+ \left(2\times \frac{1}{7}\times\frac{1}{7}\right) \right] x=12;25,25 $\ \ \ and\ \ \ $x=4;5$.

\vspace{2mm}
\noindent
(L13)\hspace{4mm}  $ \displaystyle x^2 -\left( 2+  \frac{1}{2} \right) x=13,45 $\ \ \ and\ \ \ $x=30$.

\vspace{2mm}
\noindent
(L14)\hspace{4mm}  $ \displaystyle x^2 -\left( 3+  \frac{1}{3} \right) x=13,20 $\ \ \ and\ \ \ $x=30$.   

\vspace{2mm}
\noindent
(L15)\hspace{4mm}  $ \displaystyle x^2 -\left( 4+  \frac{1}{4} \right) x=12,52;30 $\ \ \ and\ \ \ $x=30$.

\vspace{2mm}
\noindent
(L16)\hspace{4mm}  $ \displaystyle x^2 -\left( 7+  \frac{1}{7} \right) x=16,15 $\ \ \ and\ \ \ $x=35$.   

\vspace{2mm}
\noindent
(L17)\hspace{4mm}  $ \displaystyle x^2 -\left[ 7+ \left(2\times \frac{1}{7}\right) \right] x=16,10$\ \ \ and\ \ \ $x=35$.

\subsection{SMT No.\,5-Part 11-b:   $ ax-x^2  =b$}\label{SS-P11b-SMT5} 
This short subsection completes the classification of quadratic equations. For example, line 19 says

\noindent
(L19) Twice the side of a square  exceeds the area  by 0;45. (What is the side?)

\subsubsection*{Transliteration of Part 11-b}\label{SS-P11bTI-SMT5}   

\begin{note1} 
	\underline{Reverse II: Lines 18-21}\\ 
	(L18)\hspace{1mm} [1 n\'{i}gin ugu] a-\v{s}\`{a} 15 dirig\\
	(L19)\hspace{1mm} 2 n\'{i}gin ugu a-\v{s}\`{a} 45 dirig \\
	(L20)\hspace{1mm} $ \frac{\text{2}}{\text{3}} $ n\'{i}gin ugu a-[\v{s}]\`{a} 5 dirig \\
	(L21)\hspace{1mm} $ \frac{\text{1}}{\text{2}} $ n\'{i}gin \textit{ki-ma} a-[\v{s}\`{a}]-\textit{ma} 
\end{note1}

\subsubsection*{Mathematical Calculations of Part 11-b}\label{SSS-P11bMC-SMT5}

\noindent
(L18)\hspace{6mm}  $ \displaystyle  x-x^2=0;15 $\ \ \ and\ \ \ $x=0;30$.

\vspace{2mm}
\noindent
(L19)\hspace{6mm}  $ \displaystyle  2x-x^2=0;45 $\ \ \ and\ \ \ $x=0;30$.

\vspace{2mm}
\noindent
(L20)\hspace{6mm}  $ \displaystyle  \frac{2}{3} x-x^2=0;5 $\ \ \ and\ \ \ $x=0;30$.

\vspace{2mm}
\noindent
(L21)\hspace{6mm}  $ \displaystyle \frac{1}{2}x=x^2  $\ \ \ and\ \ \ $x=0;30$.

\subsection{SMT No.\,6-Part 1:    $ x^2+ax =b$}\label{SS-P1-SMT6} 

\subsubsection*{Transliteration of Part 1}\label{SSS-P1TI-SMT6}   

As a sample of the problems of this part, see the  Obverse II, lines 7-8 and   line 9 which contains the problem count.

\noindent
(L7-8) To the square of my side $11+ \frac{2}{11}$ of the length\index{length} (of my side) is added and (the result is) 1,0,40. What is the side?

\noindent
(L9) (In total) 39 addition problems (in the form of $x^2 +ax = b$).

\begin{note1} 
	\underline{Obverse I: Lines 1-30} \\
	(L1)\hspace{4mm} [a-\v{s}\`{a} nigin 1 u\v{s} dah 45] nigin \textit{mi-nu}\\
	(L2)\hspace{4mm} [a-\v{s}\`{a} nigin 2 u\v{s} dah 1],15 nigin \textit{mi-nu}\\
	(L3)\hspace{4mm} [a-\v{s}\`{a} nigin 3 u\v{s} dah]-\textit{ma} 1,45 nigin \textit{mi-nu}\\
	(L4)\hspace{4mm} [a-\v{s}\`{a} nigin 4] u\v{s} dah-\textit{ma} 2,15 nigin \textit{mi-nu}\\
	(L5)\hspace{4mm} [a-\v{s}\`{a} nigin $ \frac{\text{2}}{\text{3}} $] u\v{s} dah-\textit{ma} 35 nigin \textit{mi-nu}\\
	(L6)\hspace{4mm} [a-\v{s}\`{a} nigin $ \frac{\text{1}}{\text{2}} $] u\v{s} dah-\textit{ma} 30 nigin \textit{mi-nu}\\
	(L7)\hspace{4mm} [a-\v{s}\`{a} nigin] $ \frac{\text{1}}{\text{3}} $ u\v{s} dah-\textit{ma} 25 nigin \textit{mi-nu}\\
	(L8)\hspace{4mm} [a-\v{s}\`{a} nigin]-\textit{ia} \textbf{4} u\v{s} dah-\textit{ma} 22,30\\
	(L9)\hspace{4mm} [a-\v{s}\`{a} nigin]-\textit{ia} $ \frac{\text{1}}{\text{3}} $ \textbf{4} u\v{s} dah-\textit{ma} 17,30\\
	(L10)\hspace{2mm} [a-\v{s}\`{a} nig]in-\textit{ia} \textbf{7} u\v{s} dah-\textit{ma} 25,25\\
	(L11)\hspace{2mm} [a-\v{s}\`{a} nig]in-\textit{ia} 2 \textbf{7} u\v{s} dah-\textit{ma} 30,25\\
	(L12)\hspace{2mm} [a-\v{s}\`{a} nigin-\textit{i}]\textit{a} \textbf{7 7} u\v{s} dah-\textit{ma} 16,45,25\\
	(L13)\hspace{2mm} [a-\v{s}\`{a} nigin-\textit{i}]\textit{a} 2 \textbf{7 7} u\v{s} dah-\textit{ma} 16,50,25\\
	(L14)\hspace{2mm} [a-\v{s}\`{a} nigin \textbf{1}]\textbf{1} u\v{s} dah-\textit{ma} 55,25\\
	(L15)\hspace{2mm} [a-\v{s}\`{a} nigin] 2 \textbf{11} u\v{s} dah-\textit{ma} 1,11(sic),25\\
	(L16)\hspace{2mm} [a-\v{s}\`{a} nigin] \textbf{11 11} [u\v{s} dah-\textit{ma} 1,41,45,25]\\
	(L17)\hspace{2mm} [a-\v{s}\`{a} nigin] 2 \textbf{11 11} u\v{s} dah 1,41,50,2[5]\\
	(L18)\hspace{2mm} [a-\v{s}\`{a} nig]in \textbf{11 7} u\v{s} dah 41,15,25 \\
	(L19)\hspace{2mm} [a-\v{s}\`{a} nig]in 2 \textbf{11 7} u\v{s} dah-\textit{ma} 41,20,25\\
	(L20)\hspace{2mm} [a-\v{s}\`{a}] nigin $<$1$>$ $ \frac{\text{2}}{\text{3}} $ u\v{s} dah-\textit{ma} 1,5\\
	(L21)\hspace{2mm} [a-\v{s}\`{a} nig]in 1 $ \frac{\text{1}}{\text{2}} $ u\v{s} dah-\textit{ma} 1\\
	(L22)\hspace{2mm} [a-\v{s}\`{a} nig]in 1 $ \frac{\text{1}}{\text{3}} $ u\v{s} dah-\textit{ma} 55\\
	(L23)\hspace{2mm} [a-\v{s}\`{a} nig]in 1 \textbf{4} u\v{s} dah-\textit{ma} 52,30\\
	(L24)\hspace{2mm} [a-\v{s}\`{a} nig]in 1 $ \frac{\text{1}}{\text{3}} $ \textbf{4} u\v{s} dah-\textit{ma} 47,30\\
	(L25)\hspace{2mm} [a-\v{s}\`{a}] nigin 1 \textbf{7} u\v{s} dah 1,0,25\\
	(L26)\hspace{2mm} [a-\v{s}\`{a}] nigin 1 2 \textbf{7} u\v{s} dah 1,5,25\\
	(L27)\hspace{2mm} [a-\v{s}\`{a}] nigin 1 \textbf{7 7} u\v{s} dah-\textit{ma} 20,50,25\\
	(L28)\hspace{2mm} [a-\v{s}\`{a}] nigin 1 2 \textbf{7 7} u\v{s} dah 20,55,25\\
	(L29)\hspace{2mm} [a-\v{s}\`{a} nig]in 1 \textbf{11} u\v{s} dah [1,52,40],25\\
	(L30)\hspace{2mm} [a-\v{s}]\`{a} nigin [1 2 \textbf{11} dah 1,53,35,25]\\

	\underline{Obverse II: Lines 1-9} \\ 
	(L1)\hspace{4mm} \textit{a-na} a-\v{s}\`{a} n\'{i}gin 7 igi-7 u\v{s} dah\\
	(L2)\hspace{4mm}    24,30(sic) nigin \textit{mi}-[\textit{nu}]\\
	
	(L3)\hspace{4mm} \textit{a-na} a-\v{s}\`{a} nigin-\textit{ia} 7 2 igi-7 u\v{s}\\
	(L4)\hspace{4mm}    dah-\textit{ma} 24,35(sic) nigin \textit{mi-nu}\\
	
	(L5)\hspace{4mm} \textit{a-na} a-\v{s}\`{a} nigin-\textit{ia} 11 $<$igi-11$>$ u\v{s} dah-\textit{ma}\\
	(L6)\hspace{4mm}      1,0,35 nigin \textit{mi-nu}\\
	
	(L7)\hspace{4mm} \textit{a-na} a-\v{s}\`{a} nigin-\textit{ia} 11 2 igi-11 $<$u\v{s}$>$ dah-\textit{ma}\\
	(L8)\hspace{4mm}      1,0,40 nigin \textit{mi-nu}\\
	
	(L9)\hspace{4mm} 39 dah
	
\end{note1}

\subsubsection*{Mathematical Calculations of Part 1}\label{SSS-P1MC-SMT6}

\noindent
\underline{Obverse I, Lines 1-30:} 

\noindent
(L1)\hspace{9mm}  $ \displaystyle   x^2+x=0;45 $\ \ \ and\ \ \ $x= 0;30$.

\vspace{2mm}
\noindent
(L2)\hspace{9mm}  $ \displaystyle   x^2+2x=1;15 $\ \ \ and\ \ \ $x= 0;30$.

\vspace{2mm}
\noindent
(L3)\hspace{9mm}  $ \displaystyle   x^2+3x=1;45 $\ \ \ and\ \ \ $x= 0;30$.

\vspace{2mm}
\noindent
(L4)\hspace{9mm}  $ \displaystyle   x^2+4x=2;15 $\ \ \ and\ \ \ $x= 0;30$.

\vspace{2mm}
\noindent
(L5)\hspace{9mm}  $ \displaystyle   x^2+\frac{2}{3}x=0;35 $\ \ \ and\ \ \ $x= 0;30$.

\vspace{2mm}
\noindent
(L6)\hspace{9mm}  $ \displaystyle   x^2+\frac{1}{2}x=0;30 $\ \ \ and\ \ \ $x= 0;30$.

\vspace{2mm}
\noindent
(L7)\hspace{9mm}  $ \displaystyle   x^2+\frac{1}{3}x=0;25 $\ \ \ and\ \ \ $x= 0;30$.

\vspace{2mm}
\noindent
(L8)\hspace{9mm}  $ \displaystyle   x^2+\frac{1}{4}x=0;22,30 $\ \ \ and\ \ \ $x= 0;30$.

\vspace{2mm}
\noindent
(L9)\hspace{9mm}  $ \displaystyle   x^2+\left(\frac{1}{3}\times\frac{1}{4}\right) x=0;17,30 $\ \ \ and\ \ \ $x= 0;30$.

\vspace{2mm}
\noindent
(L10)\hspace{6mm}  $ \displaystyle   x^2+\frac{1}{7}x=0;25,25 $\ \ \ and\ \ \ $x= 0;35$.

\vspace{2mm}
\noindent
(L11)\hspace{6mm}  $ \displaystyle   x^2+\left(2\times \frac{1}{7}\right) x=0;30,25 $\ \ \ and\ \ \ $x= 0;35$.

\vspace{2mm}
\noindent
(L12)\hspace{6mm}  $ \displaystyle   x^2+\left(\frac{1}{7}\times\frac{1}{7}\right) x=16;45,25 $\ \ \ and\ \ \ $x= 4;5$.

\vspace{2mm}
\noindent
(L13)\hspace{6mm}  $ \displaystyle   x^2+\left(2\times\frac{1}{7}\times\frac{1}{7}\right) x=16;50,25 $\ \ \ and\ \ \ $x= 4;5$.

\vspace{2mm}
\noindent
(L14)\hspace{6mm}  $ \displaystyle   x^2+ \frac{1}{11}  x=0;55,25 $\ \ \ and\ \ \ $x= 0;55$.

\vspace{2mm}
\noindent
(L15)\hspace{6mm}  $ \displaystyle   x^2+\left(2\times\frac{1}{11}  \right) x=1;0,25 $\ \ \ and\ \ \ $x= 0;55$.

\vspace{2mm}
\noindent
(L16)\hspace{6mm}  $ \displaystyle   x^2+\left(\frac{1}{11}\times\frac{1}{11}\right) x=1,41;45,25 $\ \ \ and\ \ \ $x= 10;5$.

\vspace{2mm}
\noindent
(L17)\hspace{6mm}  $ \displaystyle   x^2+\left(2\times\frac{1}{11}\times\frac{1}{11}\right) x=1,41;50,25 $\ \ \ and\ \ \ $x= 10;5$.

\vspace{2mm}
\noindent
(L18)\hspace{6mm}  $ \displaystyle   x^2+\left(\frac{1}{11}\times\frac{1}{7}\right) x=41;15,25 $\ \ \ and\ \ \ $x= 6;25$.

\vspace{2mm}
\noindent
(L19)\hspace{6mm}  $ \displaystyle   x^2+\left(2\times\frac{1}{11}\times\frac{1}{7}\right) x=41;20,25 $\ \ \ and\ \ \ $x= 6;25$.

\vspace{2mm}
\noindent
(L20)\hspace{6mm}  $ \displaystyle   x^2+\left(1+\frac{2}{3}\right) x=1;5 $\ \ \ and\ \ \ $x= 0;30$.

\vspace{2mm}
\noindent
(L21)\hspace{6mm}  $ \displaystyle   x^2+\left(1+\frac{1}{2}\right) x=1  $\ \ \ and\ \ \ $x= 0;30$.

\vspace{2mm}
\noindent
(L22)\hspace{6mm}  $ \displaystyle   x^2+\left(1+\frac{1}{3}\right) x=0;55 $\ \ \ and\ \ \ $x= 0;30$.

\vspace{2mm}
\noindent
(L23)\hspace{6mm}  $ \displaystyle   x^2+\left(1+\frac{1}{4}\right) x=0;52,30 $\ \ \ and\ \ \ $x= 0;30$.

\vspace{2mm}
\noindent
(L24)\hspace{6mm}  $ \displaystyle   x^2+\left[1+\left(\frac{1}{3}\times\frac{1}{4}\right)\right] x=0;47,30 $\ \ \ and\ \ \ $x= 0;30$.

\vspace{2mm}
\noindent
(L25)\hspace{6mm}  $ \displaystyle   x^2+\left(1+\frac{1}{7}\right) x=1;0,25 $\ \ \ and\ \ \ $x= 0;35$.

\vspace{2mm}
\noindent
(L26)\hspace{6mm}  $ \displaystyle   x^2+\left[1+\left(2\times\frac{1}{7}\right)\right] x=1;5,25 $\ \ \ and\ \ \ $x= 0;35$.

\vspace{2mm}
\noindent
(L27)\hspace{6mm}  $ \displaystyle   x^2+\left[1+\left(\frac{1}{7}\times\frac{1}{7}\right)\right] x=20;50,25 $\ \ \ and\ \ \ $x= 4;5$.

\vspace{2mm}
\noindent
(L28)\hspace{6mm}  $ \displaystyle   x^2+\left[1+\left(2\times\frac{1}{7}\times\frac{1}{7}\right)\right] x=20;55,25 $\ \ \ and\ \ \ $x= 4;5$.

\vspace{2mm}
\noindent
(L29)\hspace{6mm}  $ \displaystyle   x^2+\left(1+ \frac{1}{11}\right) x=1,52;40,25 $\ \ \ and\ \ \ $x= 10;5$.

\vspace{2mm}
\noindent
(L30)\hspace{6mm}  $ \displaystyle   x^2+\left[1+ \left(2\times \frac{1}{11}\right)\right] x=1,53;35,25 $\ \ \ and\ \ \ $x= 10;5$. \\

\noindent
\underline{Obverse II, Lines 1-9:} \\ 
(L1-2)\hspace{6mm}  $ \displaystyle   x^2+\left(7+   \frac{1}{7}\right) x=24,35 $\ \ \ and\ \ \ $x= 35$.

\vspace{2mm}
\noindent
(L3-4)\hspace{6mm}  $ \displaystyle   x^2+\left[7+ \left(2\times \frac{1}{7}\right)\right] x=24,40 $\ \ \ and\ \ \ $x= 35$.

\vspace{2mm}
\noindent
(L5,6)\hspace{6mm}  $ \displaystyle   x^2+\left(11+   \frac{1}{11}\right) x=1,0,35 $\ \ \ and\ \ \ $x= 55$.

\vspace{2mm}
\noindent
(L7-8)\hspace{6mm}  $ \displaystyle   x^2+\left[11+  \left(2\times \frac{1}{11}\right)\right] x=1,0,40 $\ \ \ and\ \ \ $x= 55$.

\subsection{SMT No.\,6-Part 2:   $ x^2-ax =b$}\label{SS-P2-SMT6} 
The second part begins with the statement of a problem:\\
From the square of my side one length\index{length} (of my side) is subtracted, and (the result is) 14,30. (What is my side?)\\ 
The last line contains the problem count:\\
(In total) 30 subtraction problems (in the form of $x^2 - ax = b$).

\subsubsection*{Transliteration of Part 2}\label{SSS-P2TI-SMT6}  
\begin{note1} 
	\underline{Obverse II: Lines 10-27} \\ 
	(L10)\hspace{2mm} \textit{i-na} a-\v{s}\`{a} n\'{i}gin-\textit{ia} 1 u\v{s} zi-\textit{ma} 1[4,30]\\
	(L11)\hspace{2mm} \textit{i-na} a-\v{s}\`{a} nigin-\textit{ia} $<$2$>$ u\v{s} zi-\textit{ma} 1[4]\\
	(L12)\hspace{2mm} \textit{i-na} a-\v{s}\`{a} nigin 3 u\v{s} zi-\textit{ma} 13,30\\
	(L13)\hspace{2mm} \textit{i-na} a-\v{s}\`{a} nigin 4 u\v{s} zi-\textit{ma} 13\\
	(L14)\hspace{2mm} \textit{i-na} a-\v{s}\`{a} nigin $ \frac{\text{2}}{\text{3}} $ u\v{s} zi-\textit{ma} 14,[40]\\
	(L15)\hspace{2mm} \textit{i-na} a-\v{s}\`{a} nigin $ \frac{\text{1}}{\text{2}} $ u\v{s} zi-\textit{ma} 14,[4]5\\
	(L16)\hspace{2mm} [\textit{i}]-\textit{na} a-\v{s}\`{a} nigin $ \frac{\text{1}}{\text{3}} $ u\v{s} zi-\textit{ma} [14,50]\\
	(L17)\hspace{2mm} [\textit{i}]-\textit{na} a-\v{s}\`{a} nigin 15 u\v{s} zi-\textit{ma} 14,[5]2,30\\
	(L18)\hspace{2mm} [\textit{i-na}] a-\v{s}\`{a} nigin $ \frac{\text{1}}{\text{3}} $ \textbf{4} u\v{s} zi 14,57,30\\
	(L19)\hspace{2mm} [\textit{i-na}] a-\v{s}\`{a} nigin \textbf{7} u\v{s} zi-\textit{ma} 15,24(sic)\\
	(L20)\hspace{2mm} [\textit{i-na}] a-\v{s}\`{a} nigin-\textit{ia} 2 \textbf{7} u\v{s} zi-\textit{ma} 10,25\\
	(L21)\hspace{2mm} [\textit{i-na}] a-\v{s}\`{a} nigin \textbf{7 7} u\v{s} zi 16,35,25\\
	(L22)\hspace{2mm} \textit{i-na} a-\v{s}\`{a} nigin 2 \textbf{7 7} u\v{s} zi-[\textit{ma} 1]6,30,25\\
	(L23)\hspace{2mm} \textit{i-na} a-\v{s}\`{a} nigin \textbf{11 11} u\v{s} z[i-\textit{ma} 1,41,35,2]5\\
	(L24)\hspace{2mm} \textit{i-na} a-\v{s}\`{a} nigin 2 \textbf{11} [\textbf{1}]\textbf{1} u\v{s} [zi-\textit{ma} 1,41,30,25]\\
	(L25)\hspace{2mm} \textit{i-na} a-\v{s}[\`{a} nigin $\ldots$  $\ldots$ $\ldots$]\\
	(L26)\hspace{2mm} \textit{i-na} a-[\v{s}\`{a} nigin $\ldots$ $\ldots$ $\ldots$ $\ldots$]\\
	(L27)\hspace{2mm} \textit{i}-[\textit{na} a-\v{s}\`{a} nigin $\ldots$ $\ldots$ $\ldots$]\\

	\underline{Reverse I: Lines 1-14} \\
	(L1)\hspace{4mm} \textit{i}-[\textit{na} a-\v{s}\`{a} nigin $\ldots$ $\ldots$ $\ldots$]\\
	
	(L2)\hspace{4mm}  \textit{i-na} a-\v{s}\`{a} nigin 3 $ \frac{\text{1}}{\text{3}} $ [u\v{s} zi-\textit{ma}]\\
	(L3)\hspace{4mm}   13,20\\
	
	(L4)\hspace{4mm}  \textit{i-na} \{1\} a-\v{s}\`{a} nigin 4 u\v{s} \textbf{4} u\v{s} z[i-\textit{ma}] \\
	(L5)\hspace{4mm}   12,57(sic),30\\
	
	(L6)\hspace{4mm}  \textit{i-na} a-\v{s}\`{a} nigin-\textit{ia} 7 igi-7 u\v{s} [zi-\textit{ma}] \\
	(L7)\hspace{4mm}   16,15\\
	
	(L8)\hspace{4mm}  \textit{i-na} [a]-\v{s}\`{a} nigin-\textit{ia} 7 $<$2$>$ igi-7 $<$u\v{s}$>$ zi-\textit{ma} \\
	(L9)\hspace{4mm}   16,10\\
	
	(L10)\hspace{2mm}  \textit{i-na} a-\v{s}\`{a} nigin-\textit{ia} 11 igi-11 zi-\textit{ma}  \\
	(L11)\hspace{2mm}   40,1[5]\\
	
	(L12)\hspace{2mm}  \textit{i-na} [a-\v{s}\`{a} nigin 1]1 u\v{s} 2 ig[i-11 zi-\textit{ma}]  \\
	(L13)\hspace{2mm}   [40,10]\\
	(L14)\hspace{2mm}   30 [zi]
\end{note1}

\subsubsection*{Mathematical Calculations of Part 2}\label{SSS-P2MC-SMT6}

\noindent
\underline{Obverse II, Lines 10-27:}  

\vspace{2mm}
\noindent
(L10)\hspace{6mm}  $ \displaystyle   x^2-x=14,30 $\ \ \ and\ \ \ $x= 30$.

\vspace{2mm}
\noindent
(L11)\hspace{6mm}  $ \displaystyle   x^2-2x=14,0 $\ \ \ and\ \ \ $x= 30$.

\vspace{2mm}
\noindent
(L12)\hspace{6mm}  $ \displaystyle   x^2-3x=13,30 $\ \ \ and\ \ \ $x= 30$.

\vspace{2mm}
\noindent
(L13)\hspace{6mm}  $ \displaystyle   x^2-4x=13,0 $\ \ \ and\ \ \ $x= 30$.

\vspace{2mm}
\noindent
(L14)\hspace{6mm}  $ \displaystyle   x^2-\frac{2}{3} x=14,40 $\ \ \ and\ \ \ $x= 30$.

\vspace{2mm}
\noindent
(L15)\hspace{6mm}  $ \displaystyle   x^2-\frac{1}{2}x=14,45 $\ \ \ and\ \ \ $x= 30$.

\vspace{2mm}
\noindent
(L16)\hspace{6mm}  $ \displaystyle   x^2-\frac{1}{3}x=14,50 $\ \ \ and\ \ \ $x= 30$.

\vspace{2mm}
\noindent
(L17)\hspace{6mm}  $ \displaystyle   x^2-(0;15)x=14,52;30 $\ \ \ and\ \ \ $x= 30$.

\vspace{2mm}
\noindent
(L18)\hspace{6mm}  $ \displaystyle   x^2-\left(\frac{1}{3}\times\frac{1}{4} \right) x=14,57;30 $\ \ \ and\ \ \ $x= 30$.

\vspace{2mm}
\noindent
(L19)\hspace{6mm}  $ \displaystyle   x^2-\frac{1}{7}x=0;15,25 $\ \ \ and\ \ \ $x= 0;35$.

\vspace{2mm}
\noindent
(L20)\hspace{6mm}  $ \displaystyle   x^2-\left(2\times\frac{1}{7} \right)x=0;10,25 $\ \ \ and\ \ \ $x= 0;35$.

\vspace{2mm}
\noindent
(L21)\hspace{6mm}  $ \displaystyle   x^2-\left(\frac{1}{7}\times\frac{1}{7} \right)x=16;35,25 $\ \ \ and\ \ \ $x= 4;5$.

\vspace{2mm}
\noindent
(L22)\hspace{6mm}  $ \displaystyle   x^2-\left(2\times\frac{1}{7}\times\frac{1}{7} \right)x=16;30,25 $\ \ \ and\ \ \ $x= 4;5$.

\vspace{2mm}
\noindent
(L23)\hspace{6mm}  $ \displaystyle   x^2-\left(\frac{1}{11}\times\frac{1}{11} \right)x=1,41;35,25 $\ \ \ and\ \ \ $x= 10;5$.

\vspace{2mm}
\noindent
(L24)\hspace{6mm}  $ \displaystyle   x^2-\left(2\times\frac{1}{11}\times\frac{1}{11} \right)x=1,41;30,25 $\ \ \ and\ \ \ $x= 10;5$. \\

\noindent
\underline{Reverse I, Lines 1-14:} \\

\noindent
(L2-3)\hspace{6mm}  $ \displaystyle   x^2-\left(3+\frac{1}{3} \right)x=13,20 $\ \ \ and\ \ \ $x= 30$. \\

\vspace{2mm}
\noindent
(L4-5)\hspace{6mm}  $ \displaystyle   x^2-\left(4+\frac{1}{4} \right)x=12,52;30 $\ \ \ and\ \ \ $x= 30$. \\

\vspace{2mm}
\noindent
(L6,7)\hspace{6mm}  $ \displaystyle   x^2-\left(7+\frac{1}{7} \right)x=16,15 $\ \ \ and\ \ \ $x= 35$. \\

\vspace{2mm}
\noindent
(L8-9)\hspace{6mm}  $ \displaystyle   x^2-\left[7+\left(2\times\frac{1}{7}\right) \right]x=16,10 $\ \ \ and\ \ \ $x= 35$. \\

\vspace{2mm}
\noindent
(L10-11)\hspace{3mm}  $ \displaystyle   x^2-\left(11+\frac{1}{11} \right)x=40,15 $\ \ \ and\ \ \ $x= 55$. \\

\vspace{2mm}
\noindent
(L12-13)\hspace{3mm}  $ \displaystyle   x^2-\left[11+\left(2\times\frac{1}{11}\right) \right]x=40,10 $\ \ \ and\ \ \ $x= 55$. \\

\subsection{SMT No.\,20: First Problem}\label{SS-P1-SMT20}    
\subsubsection*{Transliteration}\label{SSS-P1TI-SMT20}  

\begin{note1} 
	\underline{Obverse:  Lines 1-12}\\
	(L1)\hspace{2mm} \textit{a-pu-s\`{a}-am-mi-ik-k}[\textit{um} a-\v{s}\`{a} \textit{\`{u}} u\v{s} ul-gar]\\
	(L2)\hspace{2mm} 36,40 [\textit{mi-nu} u\v{s}]\\
	(L3)\hspace{2mm} za-e 36,40 gar 2[6,40 igi-gub]\\
	(L4)\hspace{2mm} \textit{a-na} 36,40 \textit{i-\v{s}\'{i}-ma} [16,17,46,40 \textit{ta-mar}]\\
	(L5)\hspace{2mm} \textit{tu-\'{u}r-ma} 1 u\v{s} 1/2 \textit{h}[\textit{e-pe} 30 \textit{ta-mar}]\\
	(L6)\hspace{2mm} 30 nigin 15 \textit{ta-mar} 1[5 \textit{a-na} 16,17,46,40 dah]\\
	(L7)\hspace{2mm} 31,17,46,40 [\textit{ta-mar mi-na} \'{i}b-si]\\
	(L8)\hspace{2mm} 43,20 \'{i}b-si 30 [\textit{ta-ki-il-ta}]\\
	(L9)\hspace{2mm} \textit{i-na} 43,20 zi-\textit{ma} [13,20 \textit{ta-mar}]\\
	(L10)\hspace{0mm} igi-26,40 igi-[gub \textit{pu-\c{t}\'{u}-\'{u}r}]\\
	(L11)\hspace{0mm} 2,15 \textit{ta-mar} 2,1[5 \textit{a-na} 13,20 \textit{i-\v{s}\'{i}-ma}]\\
	(L12)\hspace{0mm} 30 \textit{ta-mar} [30 u\v{s}]
\end{note1}

\subsubsection*{Translation}\label{SSS-P1TR-SMT20}

\underline{Obverse:  Lines 1-12}
\begin{tabbing}
	\hspace{14mm} \= \kill 
	(L1)\> \tabfill{The (figure called) {\fontfamily{qpl}\selectfont \textit{apusamikkum}}. I added the area and the length (of a quadrant) together, (and the result is)}\index{length}\index{area}\index{quadrant}\index{apusamikkum@\textit{apusamikkum} (geometrical figure)}\\ 
	(L2)\> \tabfill{0;36,40. What is the length?}\index{length}\\
	(L3)\> \tabfill{You, put down 0;36,40, (and) 0;26,40 of the (area) constant.}\index{area}\\
	(L4)\> \tabfill{Multiply (0;26,40) by 0;36,40, and you see 0;16,17,46,40.}\\ 
	(L5)\> \tabfill{Return. Halve 1 of the length, (and) you see 0;30.}\index{length}\\ 
	(L6)\> \tabfill{Square 0;30, (and) you see 0;15. Add 0;15 to 0;16,17,46,40, (and)}\\
	(L7)\> \tabfill{you see 0;31,17,46,40. What is the square root?}\index{square root}\\
	(L8)\> \tabfill{0;43,20 is the square root. 0;30 (used in) completing the square,}\index{square root}\\
	(L9)\> \tabfill{subtract from 0;43,20, and you see 0;13,20.}\\ 
	(L10)\> \tabfill{Make the reciprocal of 0;26,24 of the (area) constant, (and)}\index{reciprocal of a number}\index{area} \\
	(L11)\> \tabfill{you see 2;15. Multiply 2;15 by 0;13,20, and}\\ 
	(L12)\> \tabfill{you see 0;30. 0;30 is the length.}\index{length} 
\end{tabbing}\index{completing the square method}

\subsubsection*{Mathematical Calculations}\label{SSS-P1MC-SMT20} 
\begin{figure}[H]
	\centering
	\includegraphics[scale=1]{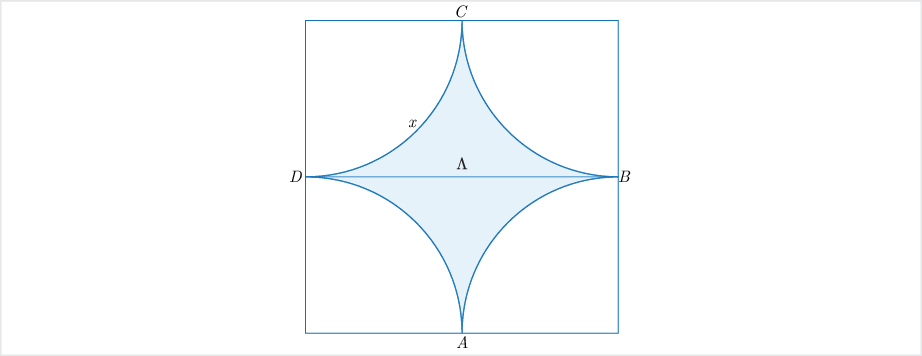}
	%\vspace{-1mm}
	\caption{Diagonal of an  {\fontfamily{qpl}\selectfont \textit{apusamikkum}}}
	\label{Figure5}
\end{figure}

Consider the {\fontfamily{qpl}\selectfont \textit{apusamikkum}}\index{apusamikkum@\textit{apusamikkum} (geometrical figure)} $ \Lambda$ as shown  in   \cref{Figure5}.   If we denote the length\index{length} of four equal quadrants of the {\fontfamily{qpl}\selectfont \textit{apusamikkum}}\index{apusamikkum@\textit{apusamikkum} (geometrical figure)} by $x$, then, as we saw in lines 22-24 of \textbf{SMT No.\,3}\index{SMT No.c@\textbf{SMT No.\,3}} (see \cite{HM22-3}),  its area\index{area of an \textit{apusamikkum}} and its diagonal\index{diagonal}  can be obtained by 
\begin{equation}\label{equ-SMT20-a}
	S_{\Lambda}\approx (0;26,40)\times x^2
\end{equation} 
and 
\begin{equation}\label{equ-SMT20-b}
	\overline{BD}\approx (1;20)\times x
\end{equation} 
respectively.  

By lines 1-2, we must have $S_{\Lambda}+\overline{BD}=0;36,40 $, so      \cref{equ-SMT20-a} and \cref{equ-SMT20-b}    give  us  the following quadratic equation\index{quadratic equation}:
\begin{equation}\label{equ-SMT20-c}
	(0;26,40) x^2 + x = 0;36,40.
\end{equation} 
The Susa scribe\index{Susa scribes} has used  completing the square\index{completing the square method} to solve this quadratic equation. According to lines 3-4, he apparently uses the new variable
\begin{equation}\label{equ-SMT20-d}
	y = (0;26,40)\times x
\end{equation}
and then multiplies both sides of \cref{equ-SMT20-c} by $0;26,40 $ to get a new quadratic equation\index{quadratic equation} with respect to new variable $y$:
\begin{equation*} 
	\left[(0;26,40)\times x\right]^2 + (0;26,40)\times x =(0;26,40)\times (0;36,40)
\end{equation*}  
or
\begin{equation}\label{equ-SMT20-e}
	y^2 + y =  0;16,17,46,40.
\end{equation} 
Now,  according to lines 3-9 
\begin{align*}
	&~~  y^2 + y =  0;16,17,46,40 \\
	\Longrightarrow~~&~~    y^2 + 2\times\left[(0;30)\times y\right] =  0;16,17,46,40 \\
	\Longrightarrow~~&~~    y^2 + 2\times\left[(0;30)\times y\right] +(0;30)^2=  0;16,17,46,40+(0;30)^2 \\
	\Longrightarrow~~&~~    y^2 + 2\times\left[(0;30)\times y\right] +(0;30)^2=  0;16,17,46,40+ 0;15  \\
	\Longrightarrow~~&~~ (y+0;30)^2  =  0;31,17,46,40 \\  
	\Longrightarrow~~&~~ y+0;30=   \sqrt{0;31,17,46,40}\\
	\Longrightarrow~~&~~ y+0;30=   \sqrt{(0;43,20)^2}\\
	\Longrightarrow~~&~~  y+0;30= 0;43,20\\
	\Longrightarrow~~&~~  y= 0;43,20-0;30
\end{align*}
which gives the solution of equation \cref{equ-SMT20-e} as
\begin{equation}\label{equ-SMT20-f}
	y=0;13,20.
\end{equation} 
Finally,  according to lines 10-12, it follows  from \cref{equ-SMT20-d} and \cref{equ-SMT20-f} that
\[x= \dfrac{1}{(0;26,40)}  \times (0;13,20)=(2;15)\times (0;13,20)=0;30 \]
which is the solution of equation \cref{equ-SMT20-c}, namely the length\index{length} of the ``side'' of the \textit{apusamikkum}\index{apusamikkum@\textit{apusamikkum} (geometrical figure)}.

\subsection{Second Problem}\label{SS-P2-SMT20}    
\subsubsection*{Transliteration}\label{SSS-P2TI-SMT20}  

\begin{note1} 
	\underline{Reverse:  Lines 13-26}\\
	(L13)\hspace{2mm} a-\v{s}\`{a}  u\v{s}  \textit{\`{u}} bar-d[\'{a} ul-gar 1,16,40]\\ 
	(L14)\hspace{2mm} za-e 26,40 igi-[gub \textit{\v{s}\`{a} ap-s\`{a}-mi-ki}]\\
	(L15)\hspace{2mm} \textit{a-na} 1,16,40 \textit{i}-[\textit{\v{s}\'{i}} 34,4,26,40 \textit{ta-mar}]\\
	(L16)\hspace{2mm} \{\textit{pu}\} \textit{tu-\'{u}r-ma} [1 u\v{s} \textit{\`{u}} 1,20 bar-d\'{a}]\\
	(L17)\hspace{2mm} \textit{\v{s}\`{a} la-a ti-du-\'{u}} [ul-gar 2,20 \textit{ta-mar}]\\
	(L18)\hspace{2mm} 1/2 2,20 \textit{he-pe} 1,[10 \textit{ta-mar} 1,10 nigin 1,21,40 \textit{ta-mar}]\\
	(L19)\hspace{2mm} 1,21,40 \textit{a-na} 34,[4,26,40 dah]\\
	(L20)\hspace{2mm} 1,55,44,[26,40 \textit{ta-mar}]\\
	(L21)\hspace{2mm} \textit{mi-na} \'{i}b-si 1,23,20 [\'{i}b-si]\\
	(L22)\hspace{2mm} 1,10 \textit{ta}-[\textit{ki-il}]-\textit{ta} \textit{i}-[\textit{na} 1,23,20 zi-\textit{ma}]\\
	(L23)\hspace{2mm} 13,2[0 \textit{ta-mar}] igi-[26,40]\\
	(L24)\hspace{2mm} \textit{pu-\c{t}\'{u}}-[\textit{\'{u}}]\textit{r} 2,15 \textit{ta}-[\textit{mar}]\\
	(L25)\hspace{2mm} 2,15 \textit{a-na} 13,20 \textit{i-\v{s}\'{i}-ma}\\
	(L26)\hspace{2mm} 30 \textit{ta-mar} 30 u\v{s}
\end{note1}

\subsubsection*{Translation}\label{SSS-P2TR-SMT20}

\underline{Reverse:  Lines 13-26}
\begin{tabbing}
	\hspace{18mm} \= \kill 
	(L13)\> \tabfill{I added the area, the length, and the diagonal, (and the result is) 1;16,40.}\index{length}\index{area}\\ 
	(L14)\> \tabfill{You, put down 0;26,40 of the (area) constant of {\fontfamily{qpl}\selectfont \textit{apusamikkum}},}\index{area of an \textit{apusamikkum}}\index{apusamikkum@\textit{apusamikkum} (geometrical figure)}\\
	(L15)\> \tabfill{multiply (it) by 1;16,40, (and) you see 0;34,4,26,49. }\\
	(L16)\> \tabfill{Return. 1 of the length and 1;20 of the diagonal}\index{length}\\ 
	(L17)\> \tabfill{that you do not know, add (them) together, (and) you see 2;20.}\\ 
	(L18)\> \tabfill{Halve 2;20, (and) you see 1;10. Square 1;10, (and) you see 1;21,40.}\\
	(L19)\> \tabfill{Add 1;21,40 to 0;34,4,26,40, (and)}\\
	(L20)\> \tabfill{you see 1;55,44,26,40.}\\
	(L21)\> \tabfill{What is the square root? 1;23,20 is the square root.}\index{square root}\\ 
	(L22)\> \tabfill{Subtract 1;10 (used in) completing the square from 1;23,20, and}\\
	(L23-24)\> \tabfill{you see 0;13,20. Make the reciprocal of 0;26,40, (and) you see 2;15.}\index{reciprocal of a number} \\ 
	(L25)\> \tabfill{Multiply 2;15 by 0;13,20, and}\\ 
	(L26)\> \tabfill{you see 0;30. 0;30 is the length.}\index{length}
\end{tabbing}\index{diagonal}\index{completing the square method} 

\subsubsection*{Mathematical Calculations}\label{SSS-P2MC-SMT20}  
If we denote the length\index{length} of four equal  quadrants of the {\fontfamily{qpl}\selectfont \textit{apusamikkum}}\index{apusamikkum@\textit{apusamikkum} (geometrical figure)} by $x$, then  line  13 and equalities  \cref{equ-SMT20-a} and \cref{equ-SMT20-b}  provide us with  the following quadratic equation\index{quadratic equation}:
\begin{equation*} 
	(0;26,40) x^2 + x+(1;20)x = 1;16,40
\end{equation*} 
or
\begin{equation}\label{equ-SMT20-g}
	(0;26,40) x^2 +  (2;20)x = 1;16,40.
\end{equation} 
Again,   according to lines 14-17, choose the new variable $y$ as  in \cref{equ-SMT20-d}
and multiply both sides of \cref{equ-SMT20-g} by $ 0;26,40$ to get a new quadratic equation\index{quadratic equation} with respect to $y$:
\begin{equation*} 
	\Big((0;26,40)\times x\Big)^2 + (2;20)\times\Big((0;26,40)\times x\Big) =(0;26,40)\times (1;16,40)
\end{equation*} 
or equivalently
\begin{equation}\label{equ-SMT20-h}
	y^2 + (2;20)y = 0;34,4,26,40.
\end{equation} 
Next,  according to lines 18-23, we write
\begin{align*}
	&~~  y^2 + (2;20)y = 0;34,4,26,40 \\
	\Longrightarrow~~&~~    y^2 + 2\times\Big((1;10)\times y\Big)=  0;34,4,26,40 \\
	\Longrightarrow~~&~~    y^2 + 2\times\Big((1;10)\times y\Big) +(1;10)^2=  0;34,4,26,40+(1;10)^2 \\
	\Longrightarrow~~&~~ (y+1;10)^2  =  0;34,4,26,40+1;21,40 \\ 
	\Longrightarrow~~&~~ (y+1;10)^2  =  1;55,44,26,40 \\   
	\Longrightarrow~~&~~ y+1;10=   \sqrt{1;55,44,26,40}\\
	\Longrightarrow~~&~~ y+1;10=   \sqrt{(1;23,20)^2}\\
	\Longrightarrow~~&~~  y+1;10= 1;23,20\\
	\Longrightarrow~~&~~  y= 1;23,20-1;10\\
	\Longrightarrow~~&~~   y= 0;13,20
\end{align*}
giving the solution of equation \cref{equ-SMT20-h}.    
Finally,   according to lines 23-26  and similar to the first problem,  \cref{equ-SMT20-d}  implies that 
\[x= \dfrac{1}{(0;26,40)}  \times (0;13,20)=(2;15)\times (0;13,20)=0;30\]
which is the solution of equation \cref{equ-SMT20-g}.

\subsection{SMT No.\,21: First Problem}\label{SS-P1-SMT21}    
\subsubsection*{Transliteration}\label{SSS-P1TI-SMT21}  

\begin{note1} 
	\underline{Obverse:  Lines 1-21}\\
	(L1)\hspace{2mm}  [\textit{a-pu-s\`{a}-mi-kum}] 5-ta-\`{a}m \textit{u\c{s}$_4$-şa-am-m}[\textit{a} nigin]
	\\
	(L2)\hspace{2mm} [\textit{ab-ni}] \textit{hi}-[\textit{in-q}]\textit{u} sag \textit{i-na li-ib-bi il}-[\textit{la-wi}]
	\\
	(L3)\hspace{2mm} [35] a-\v{s}\`{a} dal-ba-na nigin-[\textit{i}]\textit{a} \textit{mi-nu} \{\textit{nu}\}
	\\
	(L4)\hspace{2mm} [za-e] 5 \textit{mi-i}[\textit{s}]-\textit{si$_{20}$-ta a-na} 2 nigin 10 \textit{ta-mar}
	\\
	(L5)\hspace{2mm} [10] nigin 1,40 \textit{ta-mar} 1,40 \textit{i-na} 35 zi 32(sic),20 \textit{ta-mar}
	\\
	(L6)\hspace{2mm} 1 \textit{a-pu-s\`{a}-mi-ka} gar 1,20 dal \textit{\v{s}\`{a} a-pu-s\`{a}-mi-ki} gar
	\\
	(L7)\hspace{2mm} 1,20 \textit{a-na} 1 \textit{i-\v{s}\'{i}-ma} 1,20 \textit{ta-mar} 1,20 \textit{ki-ma} a-\v{s}\`{a} gar\\
	(L8)\hspace{2mm} \textit{tu-\'{u}r} 1,20 nigin 1,46,40 \textit{ta-mar} 1 nigin 1 \textit{ta-mar}
	\\
	(L9)\hspace{2mm} 1 \textit{a-na} 26,40 igi-gub \textit{a-pu-s\`{a}-mi-ki i-\v{s}\'{i}-ma}
	\\
	(L10)\hspace{0mm} 26,40 \textit{ta-mar} 26,40 \textit{i-na} 1,46,40 zi
	\\
	(L11)\hspace{0mm} 1,20 \textit{ta-mar} 1,20 \textit{a-na} 33,20 a-\v{s}\`{a} dal-ba-na \textit{i-\v{s}\'{i}-ma}
	\\
	(L12)\hspace{0mm} [4]4,26,40 \textit{ta-mar} 1,[20] \textit{a-na} 10 \textit{mi-is-si$_{20}$-ti i-\v{s}\'{i}-ma}
	\\
	(L13)\hspace{0mm} [13,20] \textit{ta-mar} 13,20 nigin [2,5]7,46,40 \textit{ta-mar}
	\\
	(L14)\hspace{0mm} [2,57,4]6,40 \textit{a-na} 4[4],26,40 dah
	\\
	(L15)\hspace{0mm} [47,24,2]6,40\textit{ t}[\textit{a-mar}] \textit{mi-na} \'{i}b-si
	\\
	(L16)\hspace{0mm} [53,20 \'{i}b-si 13,20] \textit{ta-ki-il-ta-ka}
	\\
	(L17)\hspace{0mm} [\textit{i-na} 53,20 zi 40 \textit{ta-mar} igi-1],20 \textit{\v{s}\`{a} ki-ma} a-\v{s}\`{a} gar\\
	(L18)\hspace{0mm} [\textit{pu-\c{t}\'{u}-\'{u}r} 45 \textit{ta-mar} 45 \textit{a-na} 40] \textit{i-\v{s}\'{i}-ma}
	\\
	(L19)\hspace{0mm} [30 \textit{ta-mar} u\v{s} \textit{a-pu-s\`{a}-mi-ki} 30 \textit{a-na} 1,20 dal]
	\\
	(L20)\hspace{0mm} [\textit{\v{s}\`{a} a-pu-s\`{a}-mi-ki i-\v{s}\'{i}-ma} 40 \textit{ta-mar} 40 dal]\\
	(L21)\hspace{0mm} [10 \textit{mi-is-si$_{20}$-ta a-na} 40 dah 50 \textit{ta-mar} 50 nigin]
\end{note1}

\subsubsection*{Translation}\label{SSS-P1TR-SMT21}

\underline{Obverse:  Lines 1-21}
\begin{tabbing}
	\hspace{17mm} \= \kill 
	(L1)\> \tabfill{The (figure called) {\fontfamily{qpl}\selectfont \textit{apusamikkum}}. I went (text: go) out by 5 (nindan\index{nindan (length unit)}), and I made a square (field).}\index{length}\index{apusamikkum@\textit{apusamikkum} (geometrical figure)}\\
	(L2)\> \tabfill{The narrow (that is, the {\fontfamily{qpl}\selectfont \textit{apusamikkum}}) is enclosed therein with the sides (of the inner square).}\index{apusamikkum@\textit{apusamikkum} (geometrical figure)}\\
	(L3)\> \tabfill{35,0 is the area of the space between. What is my side of the (outer) square?}\index{area of an {\fontfamily{qpl}\selectfont \textit{apusamikkum}}}\\
	(L4)\> \tabfill{You, multiply 5, the distance (between the outer and inner squares), by 2, (and) you see 10.}\\
	(L5)\> \tabfill{Square 10, (and) you see 1,40.  Subtract 1,40 from 35,0, (and) you see 33,20.}\\
	(L6)\> \tabfill{Put down 1 for the {\fontfamily{qpl}\selectfont \textit{apusamikkum}}. Put down 1;20, the diagonal of {\fontfamily{qpl}\selectfont \textit{apusamikkum}}.}\index{apusamikkum@\textit{apusamikkum} (geometrical figure)}\\
	(L7)\> \tabfill{Multiply 1;20 by 1, and you see 1;20. Put down 1;20 like the area (of the inner square).}\index{area of a square}\\
	(L8)\> \tabfill{Return. Square 1;20, (and) you see 1;46,40. Square 1, (and) you see 1.}\\ 
	(L9)\> \tabfill{Multiply 1 by 0;26,40, the constant of {\fontfamily{qpl}\selectfont \textit{apusamikkum}}, and}\\
	(L10)\> \tabfill{you see 0;26,40. Subtract 0;26,40 from 1;46,40, (and)}\\
	(L11)\> \tabfill{you see 1;20. Multiply 1;20 by 33,20, the (reduced) area of the space between, and}\index{area of an \textit{apusamikkum}}\\
	(L12)\> \tabfill{you see 44,26;40. Multiply 1;20 by 10 of the distance, and}\\
	(L13)\> \tabfill{you see 13;20. Square 13;20, (and) you see 2,57;46,40.}\\
	(L14)\> \tabfill{Add 2,57;46,40 to 44,26;40, (and)}\\
	(L15)\> \tabfill{you see 47,24;26,40. What is the square root?}\index{square root}\\
	(L16-17)\> \tabfill{53;20 is the square root. Subtract 13;20, which was (used in) your completing the square, from 53;20, (and) you see 40. The reciprocal of 1;20, 
		\{which you put down like the area,\}}\index{reciprocal of a number}\index{square root}\index{area of an \textit{apusamikkum}} \\
	(L18)\> \tabfill{make (and) you see 0;45. Multiply 0;45 by 40, and}\\
	(L19-21)\> \tabfill{you see 30. (This is) the length of {\fontfamily{qpl}\selectfont \textit{apusamikkum}}. Multiply 30 by 1;20, the diagonal of {\fontfamily{qpl}\selectfont \textit{apusamikkum}}, and you see 40. 40 is the diagonal. Add 10, the distance, to 40, (and) you see 50. 50 is the side of the (outer) square.}\index{length}
\end{tabbing}\index{diagonal}\index{completing the square method}

\subsubsection*{Mathematical Calculations}\label{SSS-P1MC-SMT21}  
\begin{figure}[H]
	\centering
	\includegraphics[scale=1]{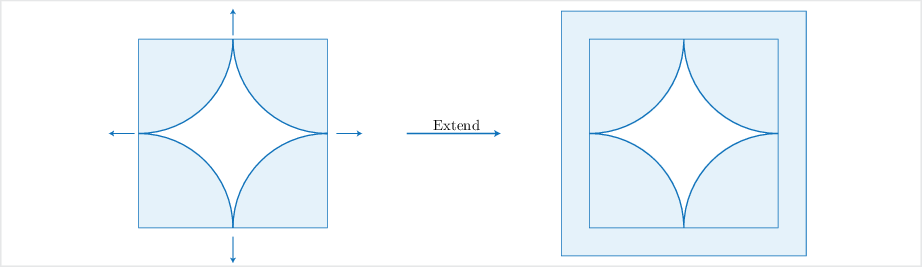}
	%\vspace{-1mm}
	\caption{Enlarging the boundary square of an {\fontfamily{qpl}\selectfont \textit{apusamikkum}}  in four directions}
	\label{Figure6}
\end{figure}

The   first problem involves the  figure {\fontfamily{qpl}\selectfont \textit{apusamikkum}}\index{apusamikkum@\textit{apusamikkum} (geometrical figure)}  already mentioned  in \textbf{SMT No.\,3} and \textbf{SMT No.\,20}\index{SMT No.t@\textbf{SMT No.\,20}} (see \cite{HM22-2} and \cref{SS-P1-SMT20}). By looking carefully at the statement of the problem in lines 1-3, we realize that   the scribe first   considers an {\fontfamily{qpl}\selectfont \textit{apusamikkum}}\index{apusamikkum@\textit{apusamikkum} (geometrical figure)}   and then enlarges each side of the boundary square\index{square} by 5 nindan\index{nindan (length unit)} in each direction  to get a bigger square\index{square} (see \cref{Figure6}). Next, he assumes that the area\index{area of a square}\index{area of an \textit{apusamikkum}} between the outer square\index{square} and the {\fontfamily{qpl}\selectfont \textit{apusamikkum}}\index{apusamikkum@\textit{apusamikkum} (geometrical figure)} to be 35,0. Finally, he asks for the value of the side of the bigger square\index{square}.

As we will see shortly, all these data result in solving a specific  quadratic equation\index{quadratic equation}.     Recall that we can construct  an  {\fontfamily{qpl}\selectfont \textit{apusamikkum}}\index{apusamikkum@\textit{apusamikkum} (geometrical figure)} by removing four equal quadrants from a square\index{square} whose side is  double the common radius of the   quadrants\index{quadrant} and whose  vertexes are     coincide with the origins of the quadrants\index{quadrant}.

\begin{figure}[H]
	\centering
	\includegraphics[scale=1]{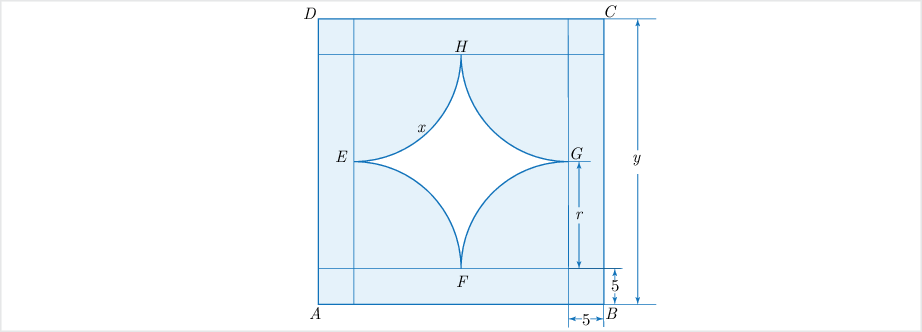}
	%\vspace{-1mm}
	\caption{The side of boundary square of an {\fontfamily{qpl}\selectfont \textit{apusamikkum}} extended in four directions}
	\label{Figure7}
\end{figure}

To formulate the problem properly, as   shown in \cref{Figure7}, we denote the  length\index{length} of each quadrant\index{quadrant} in the {\fontfamily{qpl}\selectfont \textit{apusamikkum}}\index{apusamikkum@\textit{apusamikkum} (geometrical figure)} $EFGH$ by $x$ and the length\index{length} of the bigger (or outer) square\index{square} $ABCD$ by $y$.

 It follows from the properties of  \cref{Figure7}  and the Babylonian approximation $\pi \approx  3$\index{Babylonian approximate value for $\pi$} that 
\[ 4x= 2\pi r \Longrightarrow r=\frac{4x}{2\pi}\approx \frac{4x}{6}=\frac{2x}{3} \]
and
\[ y=2(r+5)=2\times\frac{2x}{3}+10=(1;20)x+10 \]
so we get
\begin{equation}\label{equ-SMT21-a}
	\begin{dcases}
		r=\frac{2}{3}x\\
		y=(1;20)x+10.
	\end{dcases}
\end{equation}

These relations easily give us the areas of the {\fontfamily{qpl}\selectfont \textit{apusamikkum}}\index{apusamikkum@\textit{apusamikkum} (geometrical figure)} $ EFGH$ and the outer square\index{square} $ABCD $ as follows  (see \cite{HM22-3}):
\begin{equation}\label{equ-SMT21-b}
	\begin{dcases}
		S_{EFGH}=(0;26,40)x^2\\
		S_{ABCD}=\Big((1;20)x+10\Big)^2.
	\end{dcases}
\end{equation} 
Now, according to the lines 1-3, we must have
\[ S_{ABCD} - S_{EFGH}=35,0\]
which provides us with  the following quadratic equation\index{quadratic equation}
\begin{equation}\label{equ-SMT21-c}
	\Big((1;20)x+10\Big)^2 - (0;26,40)x^2=35,0.
\end{equation}  
Note that although we have obtained an equation  with respect to $x$, the scribe actually is looking for the value of the side of the outer square\index{square}, that is,  $y= (1;20)x+10$. The solving process for equation \cref{equ-SMT21-c} now begins after line 3. According to lines 4-11, we can simplify the left-hand side of \cref{equ-SMT21-c} as follows:
\begin{align*}
	&~~  \Big((1;20)x+10\Big)^2 - (0;26,40)x^2=35,0 \\
	\Longrightarrow~~&~~    \Big((1;20)x \Big)^2+2\times 10 \times  (1;20)x  +10^2 -  (0;26,40)x^2=35,0 \\
	\Longrightarrow~~&~~ (1;46,40)  x^2+2\times 10 \times  (1;20)x  +1,40 -  (0;26,40)x^2=35,0 \\
	\Longrightarrow~~&~~ (1;46,40-0;26,40)  x^2+2\times 10 \times  (1;20)x  =35,0 -    1,40 \\ 
	\Longrightarrow~~&~~  (1;20)  x^2+2\times 10 \times  (1;20)x   =33,20
\end{align*} 
so we get a new quadratic equation\index{quadratic equation}
\begin{equation}\label{equ-SMT21-d}
	(1;20)  x^2+2\times 10 \times  (1;20)x  =33,20.
\end{equation}  
From now on, the scribe  utilizes  completing the square\index{completing the square method} technique     in order to find the value of $x$ satisfying \cref{equ-SMT21-d}. According to lines 11-19, we can write
\begin{align*}
	&~~   (1;20)  x^2+2\times 10 \times  (1;20)x  =33,20 \\
	\Longrightarrow~~&~~   (1;20)\times   (1;20)  x^2 +2\times  10 \times (1;20) \times  (1;20)x   =(33,20)\times (1;20) \\
	\Longrightarrow~~&~~      \Big((1;20) x\Big)^2+2\times  (13;20) \times  (1;20)x   =44,26;40  \\
	\Longrightarrow~~&~~      \Big((1;20) x\Big)^2+2\times  (13;20) \times  (1;20)x  +(13;20)^2  =44,26;40+(13;20)^2  \\
	\Longrightarrow~~&~~      \Big((1;20) x\Big)^2+2\times  (13;20) \times  (1;20)x  +(13;20)^2  =44,26;40+2,57;46,40  \\
	\Longrightarrow~~&~~   \Big((1;20) x+13;20\Big)^2   =47,24;26,40 \\
	\Longrightarrow~~&~~    (1;20) x+13;20   = \sqrt{47,24;26,40} \\
	\Longrightarrow~~&~~    (1;20) x+13;20   = \sqrt{(53;20)^2} \\
	\Longrightarrow~~&~~    (1;20) x+13;20   = 53;20\\
	\Longrightarrow~~&~~    (1;20) x  = 53;20-13;20 \\
	\Longrightarrow~~&~~    (1;20) x  = 40 \\
	\Longrightarrow~~&~~    x  = \dfrac{1}{(1;20)} \times 40 \\
	\Longrightarrow~~&~~    x  = \left(0;45\right)\times 40
\end{align*} 
which gives us 
\begin{equation}\label{equ-SMT21-e}
	x=30.
\end{equation}  
As we said before, the scribe is interested in finding  the value of $y$. So, according to lines 19-21, we can use \cref{equ-SMT21-a} and \cref{equ-SMT21-e} to write   
\[  y=(1;20)x+10 =(1;20)\times (30)+10=40+10=50.\]
Thus, the side of the outer square\index{square} is 
$$y=50.$$

\begin{remark}\label{rem-SMT21-a}
	As is seen  in the discussion above, the scribe is extending a square\index{square} by adding to each side a value of 5 nindan\index{nindan (length unit)} in both directions. He also uses the areas of these squares\index{square} in his calculations which  require him to use and probably  expand  the algebraic expression 
	$$ \Big((1;20)x+10\Big)^2$$
	or equivalently
	$$\Big((1;20)x+2\times 5\Big)^2.$$ 
	The number $10= 2\times 5$ in this expression   suggests that the  scribes of Susa\index{scribes of Susa}  knew the algebraic  identity\index{algebraic identity}
	\begin{equation}\label{equ-SMT21-f}
		(a+2b)^2 = a^2+4ab+4b^2.
	\end{equation} 
	Because  the square of $10=2\times 5$ is   equal to  the total areas of four small squares\index{square} with sides of length\index{length} 5  in the corners of \cref{Figure7}.  In fact, a geometric explanation for verifying  the algebraic identity\index{algebraic identity} \cref{equ-SMT21-f} can be given as follows.  As is clear from \cref{Figure8},  the left-hand side of \cref{equ-SMT21-f} is the area\index{area of a square} of the outer square\index{square}  in the   figure  whose side is $a+2b$, that is $(a+2b)^2 $. 
	\begin{figure}[H]
		\centering
		\includegraphics[scale=1]{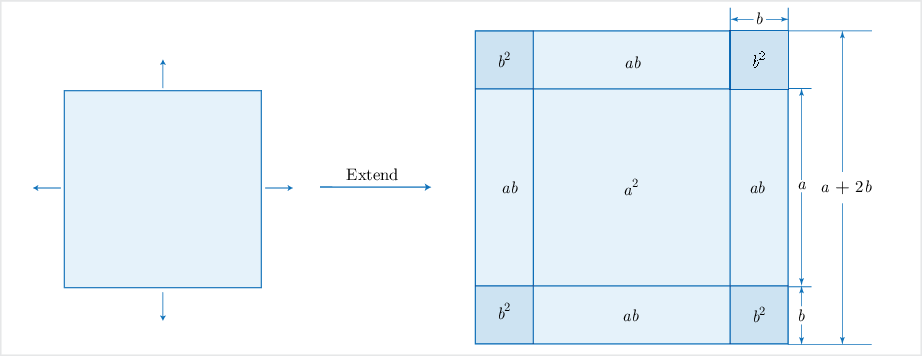}
		%\vspace{-1mm}
		\caption{Extending  a square in four directions}
		\label{Figure8}
	\end{figure}\index{square}
	
	On the other hand, the outer square\index{square} is the union of the inner square\index{square} with side $a$, four small squares\index{square} with sides $b$ and four rectangle\index{rectangle}s with  lengths $a$ and  widths $b$. Since the areas of these parts are respectively 
	\[    a^2,\  4b^2,\ \  \text{and}\  \  4ab \]
	the algebraic identity\index{algebraic identity} \cref{equ-SMT21-f} is easily verified.
\end{remark}

\subsection{SMT No.\,21: Second Problem}\label{SS-P2-SMT21}    
\subsubsection*{Transliteration}\label{SSS-P2TI-SMT21}  
\begin{note1} 
	\underline{Reverse:  Lines 1-9}\\
	(L1)\hspace{2mm}  [$ \cdots $ $ \cdots $ $ \cdots $] \textit{ta-mar}
	\\
	(L2)\hspace{2mm} [$ \cdots $ $ \cdots $ 41,40] \textit{ta-mar mi-n}[\textit{a} \'{i}]b-si\\
	(L3)\hspace{2mm} [50 \'{i}b-si 10 \textit{ta-ki-il-t}]\textit{a-ka i-na} 50 zi
	\\
	(L4)\hspace{2mm} [40 \textit{ta-mar} igi-1,20 \textit{pu}]-\textit{\c{t}\'{u}-\'{u}r} 45 \textit{ta-mar}\\
	(L5)\hspace{2mm} [45 \textit{a-na} 40 \textit{i}]-\textit{\v{s}\'{i}-ma} 30 \textit{ta-mar} u\v{s} \textit{a-pu-s\`{a}-mi-ki}\\
	(L6)\hspace{2mm} [30 \textit{a-na} 1,20] dal \textit{\v{s}\`{a} a-pu-s\`{a}-mi-ik-ki i-\v{s}\'{i}-ma}\\
	(L7)\hspace{2mm} 40 \textit{ta-mar} 40 dal 10 \textit{mi-is-si$_{20}$-ta \v{s}\`{a}} u\v{s}\\
	(L8)\hspace{2mm} \textit{a-na} 40 dah 50 \textit{ta-mar} u\v{s} 5 \textit{mi}-$<$\textit{is}$>$-\textit{si$_{20}$-ta}
	\\
	(L9)\hspace{2mm} \textit{\v{s}\`{a}} sag \textit{a-na} 40 dah 45 \textit{ta-mar} sag  
\end{note1}

\subsubsection*{Translation}\label{SSS-P2TR-SMT21}

\underline{Reverse:  Lines 1-9}
\begin{tabbing}
	\hspace{15mm} \= \kill 
	(L1)\> \tabfill{$ \cdots $ $ \cdots $ $ \cdots $ you see $ \cdots $.}\\
	(L2)\> \tabfill{$ \cdots $ $ \cdots $ $ \cdots $ you see 41,40. What is the square root?}\index{square root}\\ 
	(L3)\> \tabfill{50 is the square root. Subtract 10, which was (used in) your completing the square, from 50, (and)}\index{square root}\\
	(L4)\> \tabfill{you see 40. Make the reciprocal of 1;20, (and) you see 0;45.}\index{reciprocal of a number} \\
	(L5)\> \tabfill{Multiply 0;45 by 40, and you see 30. (This is) the length of {\fontfamily{qpl}\selectfont \textit{apusamikkum}}.}\index{length}\\
	(L6)\> \tabfill{Multiply 30 by 1;20, the diagonal of {\fontfamily{qpl}\selectfont \textit{apusamikkum}}, and}\index{apusamikkum@\textit{apusamikkum} (geometrical figure)}\\
	(L7)\> \tabfill{you see 40. 40 is the diagonal. 10, the distance of the length,}\index{length}\\
	(L8-9)\> \tabfill{add to 40, (and) you see 50. (This is) the length. 5, the distance of the width, add to 40, (and) you see 45. (This is) the width.}\index{length}\index{width}
\end{tabbing}\index{diagonal}\index{completing the square method}

\subsubsection*{Mathematical Calculations}\label{SSS-P2MC-SMT21}  
\begin{figure}[H]
	\centering
	\includegraphics[scale=1]{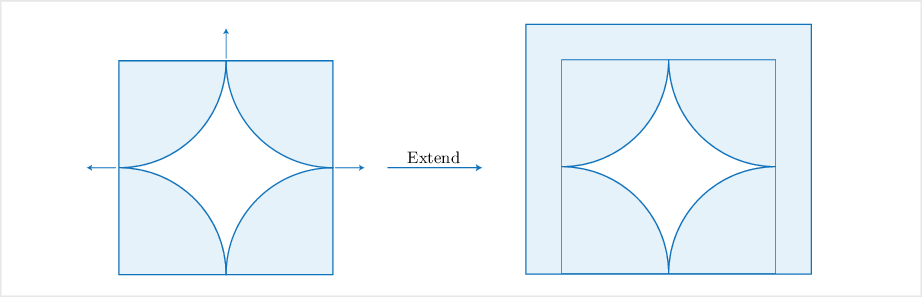}
	%\vspace{-1mm}
	\caption{Enlarging the boundary square of an {\fontfamily{qpl}\selectfont \textit{apusamikkum}}  in three directions}
	\label{Figure9}
\end{figure}

Although the beginning part of the text on the reverse is lost, we can construct the second problem and its solution by analyzing the remaining   text and comparing it with the first problem. After studying the remaining part of the second problem,  we   suggest that   the scribe  is again enlarging the boundary square\index{square} of an {\fontfamily{qpl}\selectfont \textit{apusamikkum}}\index{apusamikkum@\textit{apusamikkum} (geometrical figure)} by adding 5 nindan\index{nindan (length unit)} to the sides. Unlike the first problem, this time he appears to   extends two sides in both directions by 5 nindan\index{nindan (length unit)}  and add 5 nindan\index{nindan (length unit)} to the other two sides  in only one direction  which   produces a rectangle\index{rectangle} instead of a square\index{square} (see \cref{Figure9}). 

According to \cref{Figure9}, let us  denote the length\index{length} of each quadrant\index{quadrant}  or the side of the  our {\fontfamily{qpl}\selectfont \textit{apusamikkum}}  by $x$ and the sides of the outer rectangle\index{rectangle}  $ABCD$ by $y$ and $z$. 
Similar to the first problem, here the scribe assumes the area\index{area of a rectangle}\index{area of an \textit{apusamikkum}} between the rectangle\index{rectangle}  $ABCD$ and the  {\fontfamily{qpl}\selectfont \textit{apusamikkum}}\index{apusamikkum@\textit{apusamikkum} (geometrical figure)}  $EFGH$ to be the known value 30,50  and then asks for the sides of the rectangle\index{rectangle}.  In other words,
\[ S_{ABCD}-S_{\Lambda}=30,50. \]

\begin{figure}[H]
	\centering
	\includegraphics[scale=1]{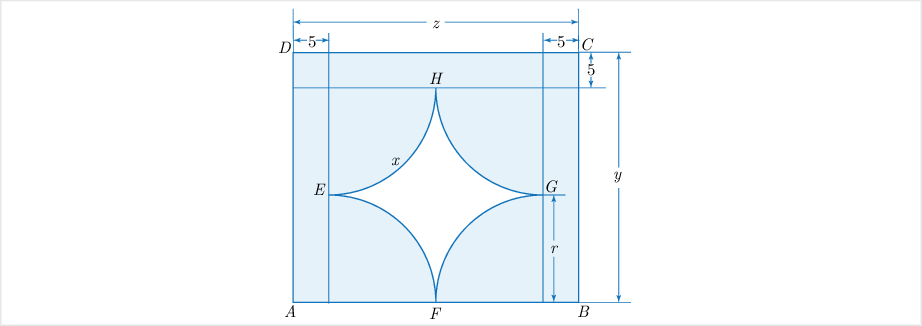}
	%\vspace{-1mm}
	\caption{The side of   boundary square of an {\fontfamily{qpl}\selectfont \textit{apusamikkum}}  extended in three directions}
	\label{Figure10}
\end{figure} 

Similar to the first problem, we have $r=(1;20)x$, so  it follows from \cref{Figure10} that 
\begin{equation}\label{equ-SMT21-h}
	\begin{dcases}
		y=(1;20)x+5\\
		z=(1;20)x+10
	\end{dcases}
\end{equation}
so the area\index{area of a rectangle}\index{area of an \textit{apusamikkum}} of the outer rectangle\index{rectangle} is 
$$ S_{ABCD}=yz= \Big((1;20)x+5\Big) \Big((1;20)x+10\Big).$$
Since the area\index{area of a rectangle}\index{area of an \textit{apusamikkum}} between the outer rectangle\index{rectangle} and the  {\fontfamily{qpl}\selectfont \textit{apusamikkum}}\index{apusamikkum@\textit{apusamikkum} (geometrical figure)}   is equal to  30,50 and $S_{EFGH}=(0;26,40)x^2$, we get the following quadratic equation\index{quadratic equation}:

\begin{equation}\label{equ-SMT21-i}
	\Big((1;20)x+5\Big) \Big((1;20)x+10\Big)-(0;26,40)x^2=30,50. 
\end{equation}
To solve equation  \cref{equ-SMT21-i}, we first need to simplify its left-hand side. So we   write 
\begin{align*}
	&~~   \Big((1;20)x+5\Big) \Big((1;20)x+10\Big)-(0;26,40)x^2=30,50 \\
	\Longrightarrow~~&~~  \Big((1;20)x\Big)^2 +10\times  [(1;20)x] + 5\times  [(1;20)x] +5\times 10-(0;26,40)x^2=30,50   \\
	\Longrightarrow~~&~~    \Big((1;20)^2 - 0;26,40\Big)    x^2 + 15\times   (1;20)x   +   50 =30,50     \\
	\Longrightarrow~~&~~   (1;46,40 - 0;26,40) x^2 +15\times  (1;20)x   =30,50-50      \\
	\Longrightarrow~~&~~    (1;20) x^2 +15\times   (1;20)x    =30,0
\end{align*}   
so we obtain the new quadratic equation\index{quadratic equation}
\begin{equation}\label{equ-SMT21-j}
	(1;20) x^2 +15\times (1;20)x   =30,0.
\end{equation} 
Similar to the first problem, the scribe  has used the usual method of completing the square\index{completing the square method} to solve this quadratic. To do so, he needs to multiply both sides of  \cref{equ-SMT21-j}  by 1;20. So we have
\begin{align*}
	&~~  (1;20) x^2 +15\times   (1;20)x   =30,0 \\
	\Longrightarrow~~&~~  (1;20)^2 x^2 +15\times (1;20) \times  (1;20)x   =(1;20)\times(30,0)   \\
	\Longrightarrow~~&~~    \Big((1;20)  x\Big)^2 +2\times (7;30) \times (1;20)  \times  (1;20)x    =40,0    \\
	\Longrightarrow~~&~~    \Big((1;20)  x\Big)^2 +2\times 10 \times (1;20)x   =40,0    \\
	\Longrightarrow~~&~~    \Big((1;20)  x\Big)^2 +2\times 10 \times  (1;20)x  +10^2  =40,0+10^2    \\
	\Longrightarrow~~&~~    \Big((1;20)  x\Big)^2 +2\times 10 \times  (1;20)x   +10^2  =40,0+1,40    \\
	\Longrightarrow~~&~~    \Big((1;20)  x+10\Big)^2   =41,40   
\end{align*}
hence we get 
\begin{equation}\label{equ-SMT21-k}
	\Big((1;20)  x+10\Big)^2   =41,40.
\end{equation} 
(Note that the value 41,40 appears in line 2 on the reverse.) Now, according to lines 3-5, the scribe finds the solution of \cref{equ-SMT21-k} as follows:
\begin{align*}
	&~~  \Big((1;20)  x+10\Big)^2   =41,40 \\
	\Longrightarrow~~&~~  (1;20)  x+10   =\sqrt{41,40}   \\
	\Longrightarrow~~&~~  (1;20)  x+10   =\sqrt{50^2}   \\
	\Longrightarrow~~&~~      (1;20)  x+10   =50   \\
	\Longrightarrow~~&~~   (1;20)  x  =50-10      \\
	\Longrightarrow~~&~~    (1;20)  x  =40    \\
	\Longrightarrow~~&~~       x  = \dfrac{1}{(1;20)} \times 40    \\ 
	\Longrightarrow~~&~~       x  = (0;45)\times 40    \\ 
	\Longrightarrow~~&~~       x  = 30. 
\end{align*}
Therefore, we get
\begin{equation}\label{equ-SMT21-l}
	x=30.
\end{equation} 
Next, according to lines 6-9, the scribe tries to find the values of $y$ and $z$. We can do this by using \cref{equ-SMT21-h} and \cref{equ-SMT21-l} as follows: 
\[ y=(1;20)x+5=(1;20)\times 30+5=40+5=45 \]
and
\[ z=(1;20)x+10=(1;20)\times 30+10=40+10=50. \]

\section{Conclusion} 
The texts of \textbf{SMT No.\,5} and \textbf{SMT No.\,6}   show  that  the    Susa scribes  of these two tablet  had good skills in working with quadratic equations. As only the statements of problems have been provided in these tablets, one may suggest that they   prepared these texts  for educational goals.  The scribes     have taken great  care   in these texts to make sure that a greater number is never subtracted from a smaller one and  to ensure that the biggest number in a sum is always mentioned first. These two facts (also observed by  Friberg  and Al-Rawi   in \cite{FA16}) imply that the concept of negative numbers was not known to the Susa scribes.  In addition, the forms of the equations in  Parts 6-12 of  \textbf{SMT No.\,5} suggest that the concept of number 0  was not known to them either. In fact,  instead of using the standard form $x^2+ax-b=0$, they have considered the alternative form $x^2+ax=b$.

Although \textbf{SMT No.\,20} and \textbf{SMT No.\,21}      deal  with quadratic equations, their problems  derive  a complicated geometric figure and the areas of its different  parts. This suggests   that  the Susa scribes not only worked with complicated geometric shapes and knew how to calculate their areas, but that  they could also combine   geometric   and  algebraic skills to formulate and solve  complicated equations. This suggests   it is reasonable to say that the scribes of Susa were using methods that later formed the basics of   ``geometric algebra''\footnote{This term refers to  the representation of algebraic concepts through geometric figures. In other words, the squares with sides  of length  $a$ can be thought of as geometric representations of $a^2$; rectangles with sides of lengths $a$ and $b$ can be interpreted as the products $ab$; and relationships among such objects can be interpreted as algebraic equations (see \cite{Kat09}).} in Greek mathematics.

\end{document}